\documentclass{amsart}
\usepackage{amssymb, amsmath, latexsym, amscd, graphicx, tabularx, color}
\usepackage[toc,page]{appendix}

\usepackage{rotating}

\usepackage[all]{xy}
\usepackage[dvipdfm]{hyperref}
\usepackage[all]{hypcap}

\newtheorem{theorem}{Theorem}[section]
\newtheorem{lemma}{Lemma}[section]

\newtheorem{definition}{Definition}[section]
\newtheorem{proposition}{Proposition}[section]
\newtheorem{remark}{Remark}[section]
\newtheorem{example}{Example}[section]

\newtheorem{atheorem}{Theorem A\hskip-0.5pt}
\newtheorem{alemma}{Lemma A\hskip-0.5pt}

\newtheorem{bproposition}{Proposition B\hskip-0.5pt}
\newtheorem{bdefinition}{Definition B\hskip-0.5pt}
\newtheorem{bexample}{Example B\hskip-0.5pt}
\newtheorem{bremark}{Remark B\hskip-0.5pt}
\newtheorem{btheorem}{Theorem B\hskip-0.5pt}
\def\pagenumber{1}

\begin{document}
\setcounter{page}{\pagenumber}
\newcommand{\T}{\mathbb{T}}
\newcommand{\R}{\mathbb{R}}
\newcommand{\Q}{\mathbb{Q}}
\newcommand{\N}{\mathbb{N}}
\newcommand{\Z}{\mathbb{Z}}
\newcommand{\tx}[1]{\quad\mbox{#1}\quad}
\parindent=0pt
\def\SRA{\hskip 2pt\hbox{$\joinrel\mathrel\circ\joinrel\to$}}
\def\tbox{\hskip 1pt\frame{\vbox{\vbox{\hbox{\boldmath$\scriptstyle\times$}}}}\hskip 2pt}
\def\circvert{\vbox{\hbox to 8.9pt{$\mid$\hskip -3.6pt $\circ$}}}
\def\IM{\hbox{\rm im}\hskip 2pt}
\def\COIM{\hbox{\rm coim}\hskip 2pt}
\def\COKER{\hbox{\rm coker}\hskip 2pt}
\def\TR{\hbox{\rm tr}\hskip 2pt}
\def\GRAD{\hbox{\rm grad}\hskip 2pt}
\def\RANK{\hbox{\rm rank}\hskip 2pt}
\def\MOD{\hbox{\rm mod}\hskip 2pt}
\def\DEN{\hbox{\rm den}\hskip 2pt}
\def\DEG{\hbox{\rm deg}\hskip 2pt}

\title[The Maslov index in PDEs geometry]{THE MASLOV INDEX IN PDEs GEOMETRY}

\author{Agostino Pr\'astaro}

\begin{abstract}
It is proved that the Maslov index naturally arises in the framework of PDEs geometry. The characterization of PDE solutions by means of Maslov index is given. With this respect, Maslov index for Lagrangian submanifolds is given on the ground of PDEs geometry. New formulas to calculate bordism groups of $(n-1)$-dimensional compact sub-manifolds bording via $n$-dimensional Lagrangian submanifolds of a fixed $2n$-dimensional symplectic manifold are obtained too. As a by-product it is given a new proof of global smooth solutions existence, defined on all $\mathbb{R}^3$, for the Navier-Stokes PDE.
Further, complementary results are given in Appendices concerning Navier-Stokes PDE and Legendrian submanifolds of contact manifolds.
\end{abstract}

\maketitle
\vspace{-.5cm}
{\footnotesize
\begin{center}
University of Rome La Sapienza, Roma Sapienza Foundation,\\ Department SBAI - Mathematics, Via A.Scarpa 16,
00161 Rome, Italy. \\
E-mail: {\tt agostino.prastaro@uniroma1.it}
\end{center}
}

\vspace{.08in} \noindent \textbf{Keywords}: Maslov index, PDE geometry; Bordism groups in PDEs; Lagrangian bordism groups; Global smooth solutions in Navier-Stokes PDEs; Legendrian submanifolds; Legendrian bordism groups.

\section{\bf Introduction}
\label{sec:1}

In 1965 V. P. Maslov introduced some integer cohomology classes useful to calculate phase shifts in semiclassical expressions for wave functions and in quantization conditions \cite{MASLOV}.\footnote{The Maslov's index is the index of a closed curve in a Lagrangian submanifold of a $2n$-dimensional symplectic space $V$, (coordinates $(x,y)$), calculated in a neighborhood of a caustic. (These are points of the Lagrangian manifold, where the projection on the $x$-plane has not constant rank $n$. Caustics are also called the projection on the $x$-plane of the set $\Sigma(V)\subset V$ of singular points of $V$, with respect to this projection.) See also \cite{KELLER}.} In the French translation, published in 1972 by the Gauthier-Villars, there is also a complementary article by V. I. Arnold, where new formulas for the calculation of these cohomology classes are given \cite{ARNOLD1,ARNOLD2, ARNOLD5}.\footnote{Further reformulations are given by L. Hormander \cite{HORMANDER}, J. Leray \cite{LERAY}, G. Lion and M. Vergue \cite{LION-VERGUE}, M. Kashiwara \cite{KASHIWARA-SCHAPIRA} and T. Thomas \cite{THOMAS}.} These studies emphasized the great importance of such invariants, hence stimulated a lot of mathematical work focused on characterization of Lagrangian Grassmannian, namely the smooth manifold of Lagrangian subspaces of a symplectic space. After the suggestion by Floer to express the spectral flow of a curve of self-adjoint operators by the Maslov index of corresponding curves of Lagrangian subspaces (1988), interesting results have been obtained relating Maslov index and spectral flow. (See, e.g., Yoside (1991), Nicolaescu (1995), Cappell, Lee and Miller (1996). )

In 1980 V. I. Arnold introduced also the notion of Lagrangian cobordism in symplectic topology \cite{ARNOLD3,ARNOLD4,ARNOLD-GIVENTAL}. This new notion has been also studied by Y. Eliashberg and M. Audin in the framework of the Algebraic Topology \cite{AUDIN,ELIASHBERG}. Next this approach has been generalized to higher order PDEs by A. Pr\'astaro \cite{PRA022}.\footnote{See also \cite{BIRAN-CORNEA} and references quoted there.}

In this paper we give a general method to recognize ``Maslov index" in the framework of the PDE geometry. Furthermore we utilize our Algebraic Topology of PDEs to calculate suitable Lagrangian bordism groups in a $2n$-dimensional symplectic manifold.

As a by-product of our geometric methods in PDEs, we get another proof of existence of global smooth solutions, defined on all $\mathbb{R}^3$, for the Navier-Stokes PDE, $(NS)$. This proof confirms one on the existence of global smooth solutions for $(NS)$, given in some our previous works \cite{PRA29,PRA33,PRA1,PRA3,PRA6,PRA14}.

Finally remark that we have written this work in an expository style, in order to be accessible at the most large audience of mathematicians and mathematical physicists.\footnote{For general complementary information on Algebraic Topology and Differential Topology, see, e.g., \cite{ALEXANDER,DUB-FOM-NOV,FERRY-RANICKI-ROSENBERG,HIRSCH2,PONT,RANICKI,SMALE1,SMALE2,SMALE3,STONG,SULLIVAN,SWITZER,TOGNOLI,WALL1,WALL2,WALL3,WARNER}.}

The main results are the following: Definition \ref{maslov-cycle-solution-pde} and Definition \ref{maslov-index-class-solution-pde} encoding Maslov cycles and Maslov indexes for solutions of PDEs that generalize usual ones. Theorem \ref{maslov-index-and-maslov-cycle-relation-solution-pde} giving a relation between Maslov cycles and Maslov indexes for solutions of PDEs. Theorem \ref{maslov-index-for-lagrangian-manifolds} recognizing Maslov index for any Lagrangian manifold, considered as solution of a suitable PDEs of first order. Theorem \ref{g-singular-lagrangian-bordism-groups} giving $G$-singular Lagrangian bordism groups, and Theorem \ref{weak-and-singular-lagrangian-bordism-groups-a} characterizing closed weak Lagrangian bordism groups. In Appendix B are reproduced similar results for Legendrian submanifolds of a contact manifold. Theorem A\ref{theorem-appendix-a} in Appendix A supports the method, given in Example \ref{navier-stokes-pdes-and-global-space-time-smooth-solutions}, to build smooth global solutions of the Navier-Stokes PDEs, defined on all $\mathbb{R}^3$.

\section{\bf Maslov index overview}\label{sec:2}

In this section we give an algebraic approach to Maslov index that is more useful to be recast in the framework of PDEs geometry. This approach essentially follows one given by V. I. Arnold \cite{ARNOLD1,ARNOLD2}, M. Kashiwara \cite{KASHIWARA-SCHAPIRA} and T. Thomas \cite{THOMAS}.

\begin{definition}\label{lagrangian-set}
Let $(V,\omega)$ be a symplectic $\mathbb{K}-$vector space over any field $\mathbb{K}$ (with characteristic $\not= 2$), where $\omega$ is a symplectic form. We denote by $L_{agr}(V,\omega)$ the {\em set of Lagrangian subspaces}, defined in {\em(\ref{set-lagrangian-subspaces})}.
\begin{equation}\label{set-lagrangian-subspaces}
  L_{agr}(V,\omega)=\left\{L<V\, |\, L=L^\bot\right\}
\end{equation}
with $E^\bot=\left\{v\in V\, |\, \omega(v,w)=0,\, \forall w\in E\right\}$.
\end{definition}

\begin{example}\label{example-lagrangian-set}
Let us consider the simplest example of $L_{agr}(V,\omega)$, with $V=\mathbb{R}^2$ and $\omega((x_1,y_1),(x_2,y_2))=x_1y_2-y_1x_2$. Then we get $L_{agr}(V,\omega)\cong G_{1,2}(\mathbb{R}^2)\cong\mathbb{R}P^1$.\footnote{We use notations and results reported in \cite{PRA000} about Grassmann manifolds.} Therefore, $L_{agr}(V,\omega)$ is a compact analytical manifold of dimension $1$. If we consider oriented Lagrangian spaces we get $L^+_{agr}(V,\omega)\cong G^+_{1,2}(\mathbb{R}^2)\cong S^1$. Since $\mathbb{R}P^1\cong S^1$, we get the commutative and exact diagram {\em(\ref{commutative-exact-diagram-lagramgian-spaces-set-in-plane})}.
\begin{equation}\label{commutative-exact-diagram-lagramgian-spaces-set-in-plane}
  \xymatrix@C=2cm{L^+_{agr}(V,\omega)\ar[d]\ar@{=}^(0.6){\backsim}_(0.6){\det^2}[r]&S^1\ar@{=}^{\wr}[d]\\
  L_{agr}(V,\omega)\ar[d]\ar@{=}^(0.6){\backsim}[r]&\mathbb{R}P^1\\
  0&&\\}
\end{equation}
In {\em(\ref{commutative-exact-diagram-lagramgian-spaces-set-in-plane})} $\det ^2$ denotes the isomorphism $L(\theta)\mapsto e^{i2\theta}$, $\theta\in[0,\pi)$.
One has the following cell decomposition into Schubert cells:
\begin{equation}\label{cell-decomposition-into-schubert-cells}
  L_{agr}(V,\omega)\cong \mathbb{R}\sqcup\{\infty\}=C_2\sqcup C_1
\end{equation}
where $C_2$ is the cell of dimension $1$ and $C_1$ is the cell of dimension $0$. This allows us to calculate the (co)homology spaces of $L_{agr}(V,\omega)$ as reported in {\em(\ref{co-homology-spaces-of-lagrangian-set-example})}.
\begin{equation}\label{co-homology-spaces-of-lagrangian-set-example}
  H^k(L_{agr}(V,\omega);\mathbb{Z}_2)\cong H_k(L_{agr}(V,\omega);\mathbb{Z}_2)\cong\bigoplus_{N_k}\mathbb{Z}_2=\left\{\begin{array}{ll}
                                                                                                                  \mathbb{Z}_2 &,\, 0\le k\le 1 \\
                                                                                                                  0&,\, k>1
                                                                                                                \end{array}\right.
                                                                                                                \end{equation}
where $N_k$ is the number of cells of dimension $k$. We get also the following fundamental homotopy group for $ L_{agr}(V,\omega)$.
\begin{equation}\label{fundamental-homotopy-group-for-example-lagranian-set }
  \pi_1(L_{agr}(V,\omega))\cong \pi_1(S^1)\cong\mathbb{Z}.
\end{equation}

$\bullet$\hskip 2pt The inverse diffeomorphism of $\det^2$, is the map $e^{i2\theta}\mapsto L(\theta)$ identifying the generator $1$ of the isomorphsism $\pi_1(L_{agr}(V,\omega))\cong\mathbb{Z}$.

$\bullet$\hskip 2pt The degree of a loop $\gamma:S^1\to L_{agr}(V,\omega)\cong S^1$, is the number of elements $\gamma^{-1}(L) $ for $L\in L_{agr}(V,\omega)$.

$\bullet$\hskip 2pt Let $\{e_1,e_2\}=\{(1,0),(0,1)\}$ be the canonical basis in $\mathbb{R}^2$. Then we call {\em real Lagrangian}
$$\mathbb{R}=\{xe_1\, |\, \forall x\in\mathbb{R}\}\subset \mathbb{R}^2$$
and {\em imaginary Lagrangian}
$$i\mathbb{R}=\{ye_2\, |\, \forall y\in\mathbb{R}\}\subset \mathbb{R}^2.$$
They are complementary: $\mathbb{R}^2\cong\mathbb{R}\bigoplus i\mathbb{R}$.

$\bullet$\hskip 2pt Let $\phi:\mathbb{R}^2\to\mathbb{R}$ be a symmetric bilinear form. One defines {\em graph} of $(\mathbb{R},\phi)$, the following set
$$\Gamma_{(\mathbb{R},\phi)}=\{(x,\phi_x(1))\in \mathbb{R}^2\}\subset\mathbb{R}^2$$
where $\phi_x:\mathbb{R}\to\mathbb{R}$ is the partial linear mapping, identified with a number via the canonical isomorphism $\mathbb{R}^*\cong\mathbb{R}$. $\Gamma_{(\mathbb{R},\phi)}$ is a Lagrangian space of $(\mathbb{R}^2,\omega)$. In fact if $x'=\lambda x$, we get $\phi_{x'}(1)=\lambda\phi_x(1)$, for any $\lambda\in\mathbb{R}$.

$\bullet$\hskip 2pt One has the identification of $L_{agr}(V,\omega)$ with a symmetric space (and Einstein manifold), via the grassmannian diffeomorphism
$$L_{agr}(V,\omega)\cong G^+_{1,2}(\mathbb{R}^2)\cong SO(2)/SO(1)\times SO(1).$$
\end{example}

\begin{example}\label{example-lagrangian-set-a}
Above considerations can be generalized to any dimension, namely considering the symeplectic space $(V,\omega)=(\mathbb{R}^{2n},\omega)$, with
$$\omega((x,y),(x',y'))=\sum _{1\le j\le n}x'_jy_j-y'_jx_j.$$

However, $L^+_{agr}(V,\omega)$ does not coincide with the grassmannian space $G^+_{n,2n}(\mathbb{R}^{2n})\cong SO(2n)/SO(n)\times SO(n)$, but one has the isomorphism reported in {\em(\ref{symmetric-space-lagrangian-example})}.\footnote{To fix ideas and nomenclature, we have reported in Tab. \ref{table-natural-geometric-structures-on-space-2n-numbers} natural geometric structures that can be recognized on $\mathbb{R}^{2n}$, besides their corresponding symmetry groups. The complex structure $i$ allows us to consider the isomorphism $\mathbb{R}^{2n}\cong\mathbb{C}^n$, $(x^j,y^j)_{1\le j\le n}\mapsto(x^j+iy^j)_{1\le j\le n}=(z^1,\cdots,z^n)$. Then the symmetry group of $(\mathbb{R}^{2n},i)\cong\mathbb{C}^n$, is $GL(n,\mathbb{C})$. Moreover the symmetry group of $(\mathbb{R}^{2n},i,\omega)$ is $Sp(n)\bigcap GL(n,\mathbb{C})=U(n)$. Therefore the matrix $A$ in (\ref{symmetric-space-lagrangian-example}) belongs to $U(n)$, hence $\det^2(A)\in\mathbb{C}$. Furthermore, taking into account that $A$ can be diagonalized with eigenvalues $\{e^{\pm i\theta_1},\cdots,e^{\pm i\theta_1}\}$, it follows that $\det^2(A)=e^{i\lambda}$ for some $\lambda\in\mathbb{R}$. Therefore, $\det^2(A)\in S^1$.}
\begin{equation}\label{symmetric-space-lagrangian-example}
  U(n)/O(n)\cong L_{agr}(V,\omega),\, A\mapsto A(i\mathbb{R}^n).
\end{equation}
Therefore one has
\begin{equation}\label{dimension-example-lagrangian-set-a}
 \dim( L_{agr}(V,\omega))=n^2-\frac{n(n-1)}{2}=\frac{n(n+1)}{2}.
\end{equation}

$\bullet$\hskip 2pt The graph $\Gamma_{(\mathbb{R}^n,\,\phi)}=\{\phi^*=\phi\in M_n(\mathbb{R})\}$ defines a chart at $\mathbb{R}^n\in L_{agr}(V,\omega)$.

$\bullet$\hskip 2pt {\em(Arnold 1967).} The square of the determinant function $\det^2:L_{agr}(V,\omega)\to S^1$, $L=A(i\mathbb{R}^n)\mapsto \det^2(A)$, induces the isomorphism
\begin{equation}\label{det-square-induced-isomorphism}
 \left\{ \begin{array}{l}
   \det{}^2_*:\pi_1(L_{agr}(V,\omega))\cong\pi_1(S^1)\cong\mathbb{Z}\\
   (\gamma:S^1\to L_{agr}(V,\omega))\mapsto {\rm degree}(\xymatrix{S^1\ar[r]^(0.4){\gamma}&L_{agr}(V,\omega)\ar[r]^(0.6){\det{}^2}&S^1\\}).
  \end{array}\right.
\end{equation}
This is a consequence of the homotopy exact sequence of the exact commutative diagram {\em(\ref{exact-commutative-diagram-fiber-bundles})} of fiber bundles.

\begin{equation}\label{exact-commutative-diagram-fiber-bundles}
  \xymatrix{&0\ar[d]&0\ar[d]&0\ar[d]&\\
  0\ar[r]&SO(n)\ar[d]\ar[r]&O(n)\ar[d]\ar[r]^{\det}&O(1)=S^0\ar[d]\ar[r]&0\\
  0\ar[r]&SU(n)\ar[d]\ar[r]&U(n)\ar[d]\ar[r]^{\det}&U(1)=S^1\ar[d]\ar[r]&0\\
  0\ar[r]&L^+_{agr}(V,\omega)\ar[d]\ar[r]&L_{agr}(V,\omega)\ar[d]\ar[r]^(0.5){\det{}^2}&L_{agr}(\mathbb{R}^2,\omega')=S^1\ar[d]\ar[r]&0\\
 &0&0&0&\\}
\end{equation}
As a by product we get the first cohomology group of $L_{agr}(V,\omega)$, with coefficients on $\mathbb{Z}$:
\begin{equation}\label{zeta-first-cohomology-group-lagrangian-set}
\left\{\begin{array}{l}
          H^1(L_{agr}(V,\omega);\mathbb{Z})=Hom_{\mathbb{Z}}(\pi_1(L_{agr}(V,\omega)),\mathbb{Z})\cong\mathbb{Z} \\
         \alpha(\gamma)={\rm degree}(\xymatrix{S^1\ar[r]^(0.4){\gamma}&L_{agr}(V,\omega)\ar[r]^(0.6){\det{}^2}&S^1\\})\in\mathbb{Z}.\\
       \end{array}\right.
\end{equation}
\end{example}

\begin{table}[t]
\caption{Natural geometric structures on $\mathbb{R}^{2n}$ and corresponding symmetry groups.}
\label{table-natural-geometric-structures-on-space-2n-numbers}
\scalebox{0.6}{$\begin{tabular}{|l|l|l|}
\hline
\hfil{\rm{\footnotesize Name}}\hfil&\hfil{\rm{\footnotesize Structure}}&\hfil{\rm{\footnotesize Symmetry group}}\hfil\\
\hline\hline
\hfil{\rm{\footnotesize Euclidean}}\hfil&\hfil{\rm{\footnotesize $g:\mathbb{R}^{2n}\times\mathbb{R}^{2n}\to\mathbb{R}$}}&\hfil{\rm{\footnotesize $O(2n)=\{A=(a_{ij})\in M_{2n}(\mathbb{R})\, |\, \det A\not =0,\, A^*A=I_{2n}\}$}}\hfil\\
\hfil{\rm{\footnotesize }}\hfil&\hfil{\rm{\footnotesize $g(v,v')=\sum_{1\le j\le n}(x_jx'_j+y_jy'_j)$}}&\hfil{\rm{\footnotesize }}\hfil\\
\hline
\hfil{\rm{\footnotesize symplectic}}\hfil&\hfil{\rm{\footnotesize $\omega:\mathbb{R}^{2n}\times\mathbb{R}^{2n}\to\mathbb{R}$}}&\hfil{\rm{\footnotesize $Sp(n)=\{A=(a_{ij})\in M_{2n}(\mathbb{R})\, |\, \det A\not =0,\, A^*\left(
                                                                      \begin{array}{cc}
                                                                        0 & I_n \\
                                                                        I_n & 0 \\
                                                                      \end{array}
                                                                    \right)
A=\left(
    \begin{array}{cc}
      0 & I_n \\
      -I_n & 0 \\
    \end{array}
  \right)
\}$}}\hfil\\
\hfil{\rm{\footnotesize }}\hfil&\hfil{\rm{\footnotesize $\omega(v,v')=\sum_{1\le j\le n}(x'_jy_j-x_jy'_j)$}}&\hfil{\rm{\footnotesize }}\hfil\\
\hline
\hfil{\rm{\footnotesize hermitian}}\hfil&\hfil{\rm{\footnotesize $h:\mathbb{R}^{2n}\times\mathbb{R}^{2n}\to\mathbb{C}$}}&\hfil{\rm{\footnotesize $U(n)=\{A=(a_{ij})\in M_n(\mathbb{C})\, |\, \det A\not =0,\, AA^*=I_n\}$}}\hfil\\
\hfil{\rm{\footnotesize }}\hfil&\hfil{\rm{\footnotesize $h(v,v')=g(v,v')+i\omega(v,v')=\sum_{1\le j\le n}(x_j+iy_j)(x'_j+iy'_j)$}}&\hfil{\rm{\footnotesize }}\hfil\\
\hline
\multicolumn{3}{l}{\rm{\footnotesize $U(n)=Sp(n)\bigcap O(2n)$. $O(2n),\, Sp(n)\subset GL(2n,\mathbb{R})$, closed sub-groups.}}\hfill\\
\multicolumn{3}{l}{\rm{\footnotesize $GL(n,\mathbb{C})$ is the symmetry group of the complex structure. $O(2n)\bigcap GL(n,\mathbb{C})=GL(n,\mathbb{C})\bigcap Sp(n)=U(n)$.}}\hfill\\
\multicolumn{3}{l}{\rm{\footnotesize $\mathbb{R}^{2n}=\{(x,y)=(x_1,\cdots,x_n,y_1,\cdots, y_n)\, |\, x_j,\, y_j\in\mathbb{R}\}$.}}\hfill\\
\multicolumn{3}{l}{\rm{\footnotesize $A^*=(\bar a_{ji})$, if $A=(a_{ij})$. In the real case $A^*=(a_{ji})$.}}\hfill\\
\end{tabular}$}
\end{table}

\begin{example}\label{example-lagrangian-set-b}
Let $(V,\sigma)$ be a $2n$-dimensional real symplectic vector space, endowed with a complex structure $J:V\to V$, such that $g:V\times V\to \mathbb{R}$, $g(u,v)=\sigma(J(u),v)$, is an inner product. Then for any $L\in L_{agr}(V,\sigma)$, the following propositions hold.

{\em(i)} One has the diffeomorphism $U(V)/O(L)\cong L_{agr}(V,\sigma)$, $A\mapsto A(L)$, where
$$O(L)=\{A\in U(V)\, |\, A(L)=L\}.$$

{\em(ii)} One has the isomorphism $f_L:(\mathbb{R}^{2n},\omega,i)\cong(V,\sigma,J)$, $f_L(i\mathbb{R}^n)=L$.

{\em(iii)} One has the diffeomorphism $$f_L:L_{agr}(\mathbb{R}^{2n},\omega)\cong U(n)/O(n)\to L_{agr}(V,\sigma)\cong U(V)/O(L),\, \lambda\mapsto f_L(\lambda).$$

{\em(iv)} If $L_1,\, L_2\in L_{agr}(V,\sigma)$, there exists a {\em difference element} $\lambda[L_1,L_2]\in L_{agr}(\mathbb{R}^{2n},\omega)$, such that $\lambda[L_1,L_2]\cong i\mathbb{R}^n\subset\mathbb{R}^{2n}$. Therefore $L_{agr}(V,\sigma)$ has a $L_{agr}(\mathbb{R}^{2n},\omega)$-affine structure.
\end{example}

\begin{definition}\label{witt-group}
The {\em Witt group} of a field $\mathbb{K}$ is $W(\mathbb{K})=\pi_0(\mathcal{Q}_+)$, where $\mathcal{Q}_+$ is the category whose objects are {\em quadratic spaces}, namely $\mathbb{K}-$vector spaces with non-degenerate, symmetric bilinear forms.
We say that two quadratic spaces $V_1,\, V_2\in Ob(\mathcal{Q}_+)$, are {\it Witt-equivalent} if there exists a Lagrangian correspondence between them, more precisely a morphism $f\in Hom_{\mathcal{Q}_+}(V_1,V_2):=L_{agr}(V_1^o\oplus V_2)$, called the {\it space of Lagrangian correspondences}. There $(V,q)^o:=(V,-q)$, with $q$ the quadratic structure. Composition of morphisms is meant in the sense of composition of general correspondences. (For example if $f:V_1\to V_2$ is an isometry then the graph $\Gamma_f\subset V_1^o\oplus V_2$ is Lagrangian. Think of composing functions $f:A\to B$ and $g:B\to C$ via the subsets of $A\times B$ and $B\times C$.)\footnote{If $f: V_1\to V_2$ is an isomorphism, then the graph $\Gamma_f\subset V_1^o\bigoplus V_2$ is Lagrangian. The quadratic space $(V,q)$ is equivalent to $0$ iff it contains Lagrangian. (For more details on Witt group see the following link: \href{http://en.wikipedia.org/wiki/Witt_group}{Wikipedia-Witt-group} and References therein.)}
\end{definition}

\begin{proposition}\label{property-witt-group}
$W(\mathbb{K})$ is the group whose elements are Witt equivalence classes of quadratic spaces, with addition induced by direct sum, and the inverse $-(V,q)$ is given by $-(V,q)=(V,q)^o$.
\end{proposition}

\begin{example}\label{example-real-witt-group}
Let us consider $$(V,q)=\left(\mathbb{K}^2,\, \left(
                                            \begin{array}{cc}
                                              +1 & 0 \\
                                              0 & -1 \\
                                            \end{array}
                                          \right)\right).$$

 $\bullet$\hskip 2pt One has the isomorphism: $W(\mathbb{R})\cong \mathbb{Z}$ that is the index of $q$, namely the number of positive eigenvalues minus the number of negative eigenvalues.

 $\bullet$\hskip 2pt One has the isomorphism: $W(\mathbb{C})\cong \mathbb{Z}/2\mathbb{Z}=\mathbb{Z}_2$ that is the dimension of $W(\mathbb{C})$.\footnote{In this paper we denote $\mathbb{Z}/n\mathbb{Z}$ by $\mathbb{Z}_n$.}
\end{example}

\begin{theorem}\label{maslov-index-of-n-tuples-lagrangian}
There exists a canonical mapping $\tau: L_{agr}(V,\omega)^{\mathbb{Z}_r}\to W(\mathbb{K})$ that we call {\em Maslov index} and that factorizes as reported in the commutative diagram {\em(\ref{commutative-diagram-factors-maslov-index})}.
\begin{equation}\label{commutative-diagram-factors-maslov-index}
  \xymatrix{L_{agr}(V,\omega)^{\mathbb{Z}_r}\ar[dr]_{\tau}\ar[r]&Ob(\mathcal{Q}_+)\ar[d]\\
  &W(\mathbb{K})\\}
\end{equation}
\end{theorem}
\begin{proof}
Given a $r$-tuple $L=(L_1,\cdots,L_r)$ of Lagrangian subspaces of $(V,\omega)$, we can identify a cochain complex (\ref{fundamental-chain-complex})
\begin{equation}\label{fundamental-chain-complex}
 \xymatrix{\framebox{$C_L=\bigoplus(L_i\bigcap L_{i+1})$}\ar[r]^(0.7){\partial}&\bigoplus_iL_i\ar[r]^{\Sigma} &V\\}
\end{equation}

where $\Sigma$ is the sum of the components, and $\partial(a)=(a,-a)\in L_i\oplus L_{i+1}$, $\forall a\in L_i\bigcap L_{i+1}$. Then we get a quadratic space $(T_L,q_L)$, with  $T_L=\ker\sum/\hbox{\rm im }\partial$ and $q_L(a,b)=\sum_{i>j}\omega(a_i,b_j)$, ({\it Maslov form}), where $a,\, b\in T_L$ are lifted to the representative $(a_i),\, (b_i)\in \oplus_{i} L_i$. Then the {\em Maslov index} is defined by  (\ref{definition-maslov-index}).

\begin{equation}\label{definition-maslov-index}
  \tau(L)=\tau(L_1,\cdots,L_r)=(T_L,q_L)\in W(\mathbb{K}).
\end{equation}
  One has the following properties:

(a) {\em Isometries}: $T(L_1,\cdots,L_r)=T(L_r,L_1,\cdots,L_{r-1})=T(L_1,\cdots,L_r)^o$.

(b) {\em Lagrangian correspondences}:

$T(L_1,\cdots,L_r)\oplus T(L_1,L_k,\cdots,L_{r})\to T(L_1,\cdots,L_r)$, $k<r$.

By considering cell complex $C_L=C(L_1,\cdots,L_r)$, as {\it $r$-gon}, with the face labelled by $V$, edges labelled by $L_i$ and vertices labelled by $L_i\bigcap L_{i+1}$, property (b) allows us to reduce to the case of three Lagrangian subspaces. Furthermore, Lagrangian correspondences induce cobordism properties. For example $C(L_1,L_2,L_3,L_4)$ cobords with $C(L_1,L_2,L_3)\bigcup C(L_1,L_3,L_4)$.

(c) {\em Cocycle property}:

$\tau(L_1,L_2,L_3)-\tau(L_1,L_2,L_4)+\tau(L_1,L_3,L_4)-\tau(L_2,L_3,L_4)=0$.

\end{proof}

\begin{theorem}[Leray's function]\label{topology-lagrangian-grassmannian}

$\bullet$\hskip 2pt {\em (Case $\mathbb{K}=\mathbb{R}$).}

Let $\pi:\widetilde{L_{agr}(V,\omega)}\to  L_{agr}(V,\omega)$ be the universal cover of the Lagrangian Grassmannian. Then there exists a function {\em( Leray's function)}
$$m:\widetilde{L_{agr}(V,\omega)}^2\to\mathbb{Z}\cong W(\mathbb{R})\cong\pi_1(L_{agr}(V,\omega)) $$
such that
$$\tau(\pi(\widetilde{L_1}),\cdots,\pi(\widetilde{L_r}))=\sum_{i\in\mathbb{Z}_r}m(\widetilde{L_i},\widetilde{L}_{i+1}).$$

$\bullet$\hskip 2pt {\em (Case $\mathbb{K}$ general ground field).}

Let $L^+_{agr}(V,\omega)$ be the set of oriented Lagrangians. There exists a function $$m:L^+_{agr}(V,\omega)\to W(\mathbb{K})$$ such that
$$\tau(L_1,\cdots,L_r)=\sum_i m(L_i,L_{i+1})\hskip 3pt \hbox{\rm mod}\hskip 3pt I^2$$
where $I=\ker(\dim:W(\mathbb{K})\to \mathbb{Z}_2)$. \footnote{$L_{agr}(V,\omega)$ has a unique double cover $L^{(2)}agr(V,\omega)$. For any pair $(\widetilde{L_1},\widetilde{L_2})$ with $\widetilde{L_1},\, \widetilde{L_2}\in L^{(2)}agr(V,\omega)$, the number $m(\widetilde{L_1},\widetilde{L_2})$ is well-defined mod $4$.}
\end{theorem}

\begin{theorem}[Metaplectic group]\label{metaplectic-group}
The Maslov index allows to identify a central extension $Mp(V)$ of the group $Sp(V)$ that when $\mathbb{K}=\mathbb{R}$ is the unique double cover of $Sp(V)$. ($Mp(V)$ is called {\em metaplectic group}.)
\end{theorem}
\begin{proof}
The cocycle property allows to equip $Mp_1(V)=W(\mathbb{K})\times Sp(V)$, with the multiplication:
\begin{equation}\label{multiplication}
  (q,g).(q',g')=(q+q'+\tau(L,gL,gg'L),gg').
\end{equation}
Thus $Mp_1(V)$ is a group and gives a central extension
\begin{equation}\label{central-extension-a}
  \xymatrix{0\ar[r]&W(\mathbb{K})\ar[r]&Mp_1(V)\ar[r]&Sp(V)\ar[r]&1\\}
\end{equation}
Moreover set
\begin{equation}\label{mp2}
 Mp_2(V)=\left\{(m(g\widetilde{L},\widetilde{L})+q,g)\, |\, q\in I^2,\, g\in Sp(V)\right\}\subset Mp_1(V)
\end{equation}
where $\widetilde{L}\in\Lambda$ over $L\in L_{agr}(V)$. $Mp_2(V)$ is a subgroup, giving a central extension
\begin{equation}\label{central-extension-b}
  \xymatrix{0\ar[r]&I^2\ar[r]&Mp_2(V)\ar[r]&Sp(V)\ar[r]&1\\}
\end{equation}
By quotient $I^2$ by $I^3$ we define a central extension
\begin{equation}\label{central-extension-c}
  \xymatrix{0\ar[r]&I^2/I^3\ar[r]&Mp(V)\ar[r]&Sp(V)\ar[r]&1\\}
\end{equation}
defining $Mp(V)$, called {\em metaplectic group}.

When $\mathbb{K}=\mathbb{R}$, $I^2/I^3\cong \mathbb{Z}_2$, so $Mp(V)$ is the unique double cover of $Sp(V)$. In this case $Mp(V)$ has four connected components, among which $Mp_2(V)$ is the identity. $Mp_2(V)$ is the universal covering group of $Sp(V)$.\footnote{One can construct $Mp(V)$ also by observing that $Sp(V)$ embeds into $L_{agr}(V^o\bigoplus V)$ by $g\mapsto \Gamma_g$, the graph of $g$. Then define multiplication on $Mp_2(V)$: $(q,g).(q',g')=(q+q'+\tau(\Gamma_1,\Gamma_g,\Gamma_{gg'}),gg')$. Moreover, $\Gamma_g$ has a canonical orientation.}
\end{proof}

\begin{example}[Arnold's Maslov index]\label{arnold-maslov-index}
The cohomology class of the Arnold's approach for Maslov index is $\alpha\in H^1(L_{agr}(\mathbb{R}^{2n},\omega);\mathbb{Z})\cong \mathbb{Z}$, obtained as the pullback of the standard differential form $d\theta:S^1\to T^*S^1$, via $\det{}^2:L_{agr}(\mathbb{R}^{2n},\omega)\to S^1$. In  {\em(\ref{arnold-definitions-maslov-index})} are summarized the Arnold's definitions of Maslov index for $L\in L_{agr}(\mathbb{R}^{2},\omega)$.\footnote{In particular, if $0\le\theta_1<\theta_2<\theta_3<\pi$, then $\tau(L_1,L_2,L_3)=1$.}

\begin{equation}\label{arnold-definitions-maslov-index}
\left\{\begin{array}{ll}
  \tau(L(\theta))& = \left\{\begin{array}{ll}
                             1-\frac{2\theta}{\pi}&,\, 0< \theta<\pi\\
                              0& ,\,\theta=0\\
                            \end{array}\right.\\
  \tau(L_1,L_2)&=- \tau(L_2,L_1)= \left\{\begin{array}{ll}
                             1-\frac{2(\theta_1-\theta2)}{\pi}&,\, 0\le \theta_1<\theta_2<\pi\\
                              0& ,\,\theta_1=\theta_2\\
                            \end{array}\right.\\
   \tau(L_1,L_2,L_3)&=\tau(L_1,L_2)+\tau(L_2,L_3)+\tau(L_3,L_1)\in\left\{-1,0,1\right\}\subset\mathbb{Z}.\\
  \end{array}\right.
\end{equation}

$\bullet$\hskip 2pt Any couple $(L_1,L_2)$ of Lagrangians in $L_{agr}(\mathbb{R}^2,\omega)$, determines a curve $\gamma_{12}:I=[0,1]\to L_{agr}(\mathbb{R}^2,\omega)$, $\gamma_{12}(t)=L((1-t)\theta_1+t\theta_2)$, connecting $L_1$ and $L_2$.

$\bullet$\hskip 2pt A triple $(L_1,L_2,L_3)$ of Lagrangians in $L_{agr}(\mathbb{R}^2,\omega)$, determines a loop $\gamma_{123}=\gamma_{12}\gamma_{23}\gamma_{31}:S^1\to L_{agr}(\mathbb{R}^2,\omega)$, with homotopy class the Maslov index of the triple:
$$\gamma_{123}=\tau(L_1,L_2,L_3)\in\{-1,0,1\}\subset\pi_1(L_{agr}(\mathbb{R}^2,\omega))\cong\mathbb{Z}.$$
In fact, for $0\le\theta_1<\theta_2<\theta_3<\pi$, one has $\det{}^2\gamma_{123}=1:S^1\to S^1$, and ${\rm degree}(\det{}^2\gamma_{123})=1=\tau(L_1,L_2,L_3)\in\mathbb{Z}$.

In  {\em(\ref{arnold-definitions-maslov-index-a})} are summarized the Arnold's definitions of Maslov index for $L\in L_{agr}(\mathbb{R}^{2n},\omega)$, $n>1$. There $\pm e^{i\theta_{1}},\cdots,\pm e^{i\theta_{n}}$, denote the eigenvalues of the matrix $A\in U(n)$, such that $A(i\mathbb{R}^n)=L$.

\begin{equation}\label{arnold-definitions-maslov-index-a}
\left\{\begin{array}{ll}
  \tau(L)& =\sum_{1\le j\le n}(1-\frac{2\theta_j}{\pi})\in \mathbb{R},\, 0\le \theta_j<\pi\\
  \tau(L_1,L_2)&=- \tau(L_2,L_1)= \left\{\begin{array}{ll}
                             \sum_{1\le j\le n}(1-\frac{2(\theta_{1j}-\theta_{2j})}{\pi}&,\, 0\le \theta_{1j}<\theta_{2j}<\pi\\
                              0& ,\,\theta_{1j}=\theta_{2j}\\
                            \end{array}\right.\\
   \tau(L_1,L_2,L_3)&=\tau(L_1,L_2)+\tau(L_2,L_3)+\tau(L_3,L_1)\in\left\{-1,0,1\right\}\subset\mathbb{Z}.\\
  \end{array}\right.
\end{equation}

$\bullet$\hskip 2pt {\em(Arnold 1967).} The Poincar\'e dual $D\alpha$ of $\alpha\in H^1(L_{agr}(\mathbb{R}^{2n},\omega);\mathbb{Z})\cong \mathbb{Z}$, is called the {\em Maslov cycle}, and it results
 \begin{equation}\label{maslov-cycle}
  D\alpha=\{L\in L_{agr}(\mathbb{R}^{2n},\omega)\, | \, L\bigcap i\mathbb{R}^n\not=\{0\}\}
\end{equation}
with
\begin{equation}\label{maslov-cycle-a}
  [D\alpha]\in H_{\frac{(n+2)(n-1)}{2}}(L_{agr}(\mathbb{R}^{2n},\omega);\mathbb{Z}).
\end{equation}
\end{example}

\begin{example}[The Wall non-additivity invariant as Maslov index]\label{relation-maslov-index-wall-non-addititivity-invariant}
Let $(V,\omega)$ be a symplectic space and $(L_1,L_2,L_3)$ a triple of Lagrangian subspaces. The {\em Wall non-additivity invariant} $w(L_1,L_2,L_3)=\sigma(W,\psi)$, i.e., the signature of the non-singular symmetric form
$$\psi:W\times W\to \mathbb{R},\, \pi(x_1,x_2,x_3,y_1,y_2,y_3)=\omega(x_1,y_2)$$
with
$$W=\frac{\{(x_1,x_2,x_3)\in L_1\oplus L_2\oplus L_3\, |\, x_1+x_2+x_3=0\}}{{\rm im}(L_1\bigcap L_2+L_2\bigcap L_3+L_3\bigcap L_1)}.$$

$\bullet$\hskip 2pt {\em(Wall \cite{WALL3})}  $w(L_1,L_2,L_3)$ can be identified with the defect of the Novikov additivity for the signature of the triple union of a $4k$-dimenaional manifold with boundary $(X,\partial X)$:
$$w(L_1,L_2,L_3)=\sigma(X_1)+\sigma(X_1)+\sigma(X_2)+\sigma(X_3)-\sigma(X)\in\mathbb{Z}$$

where $X=X_1\bigcup X_2\bigcup X_3$, and $X_i$, $i=1,2,3$, are codimension $0$ manifolds with boundary meeting transversely as pictured in {\em(\ref{relation-wall-invariant-novikov-defect})}. One has a nonsingular symplectic intersection form on $H^{2k-1}(X_1\bigcap X_2\bigcap X_3;\mathbb{R})$, \footnote{The {\em intersection form} of a $2n$-dimensional topological manifold with boundary $(M,\partial M)$, is $(-1)^{n}$-symmetric form $\lambda:H^n(M,\partial M;\mathbb{Z})/Tor\times  H^n(M,\partial M;\mathbb{Z})/Tor\to\mathbb{Z}$, $\lambda(x,y)=<x\bigcup y,[M]>\in\mathbb{Z}$. The {\em signature} $\sigma(M)$ of $4k$-dimensional manifold $(M,\partial M)$, is $\sigma(M)=\sigma(\lambda)\in \mathbb{Z}$.} and the following Lagrangian subspaces:
\begin{equation}\label{Lagrangian-subspaces-wall-invariant}
  \left\{
\begin{array}{l}
   L_1={\rm im}(H^{2k-1}(X_2\bigcap X_3;\mathbb{R})\to H^{2k-1}(X_1\bigcap X_2\bigcap X_3;\mathbb{R}))\\
     L_2={\rm im}(H^{2k-1}(X_1\bigcap X_3;\mathbb{R})\to H^{2k-1}(X_1\bigcap X_2\bigcap X_3;\mathbb{R}))\\
      L_3={\rm im}(H^{2k-1}(X_1\bigcap X_2;\mathbb{R})\to H^{2k-1}(X_1\bigcap X_2\bigcap X_3;\mathbb{R})).\\
  \end{array}
  \right.
\end{equation}

$\bullet$\hskip 2pt {\em(Cappell, Lee and Miller \cite{CAPPELL-LEE-MILLER})} The Maslov index of the triple $(L_1,L_2,L_3)$ coincides with the Wall non-dditivity invariant of $(L_1,L_2,L_3)$.\footnote{A more recent different proof has been given by A. Ranicki (1997). (See in \cite{RANICKI}.)}
$$\tau(L_1,L_2,L_3)=w(L_1, L_2,L_3,g).$$
\end{example}
\begin{equation}\label{relation-wall-invariant-novikov-defect}
 \xymatrix{&&\bullet\ar@/_3pc/@{-}[dddll]_{X_1\bigcap\partial X}\ar@{-}[d]\ar@/^3pc/@{-}[dddrr]^{X_2\bigcap\partial X}&&\\
  &X_1&\bullet\ar@{-}[ddll]\ar@{-}[ddrr]&X_2&\\
  &&X_3&&\\
  \bullet\ar@/_1pc/@{-}[rrrr]_{X_3\bigcap \partial X}&&&&\bullet\\}
\end{equation}

\section{\bf Integral bordism groups in PDEs}\label{sec:3}

The definition of Maslov index can be recast in the framework of the PDE's geometry. In fact the {\it metasymplectic structure} of the Cartan distribution of $k$-jet-spaces $J^k_n(W)$ over a fiber bundle $\pi:W\to M$, $\dim W=n+m$, $\dim M=n$, allows us to recognize ``Maslov index" associated to $n$-dimensional integral planes of the Cartan distribution of $J^k_n(W)$, and by restriction on any PDE $E_k\subset J^k_n(W)$. In the following we shall give a short panorama on the geometric theory of PDEs and on the metasymplectic structure of the Cartan distribution and its relations with (singular) solutions of PDEs. (For more information see also \cite{LYCH-PRAS,PRA1}.)\footnote{For general information on PDE's geometry see \cite{B-C-G-G-G,HCARTAN,GOLD1,GROMOV,KRAS-LYCH-VIN,PRA000}.}

Let $W$ be a smooth manifold of dimension $m+n$. For any
$n$-dimensional submanifold $N\subset W$ we denote by $[N]^k_a$ the
{\em $k$-jet of $N$ at the point $a\in N$}, i.e., the set of $n$-dimensional
submanifolds of $W$ that have in $a$ a contact of order $k$. Set
$J^k_n(W)\equiv\bigcup_{a\in W}J^k_n(W)_a$, $J^k_n(W)_a\equiv\{[N]^k_a | a\in W\}$.
We call $J^k_n(W)$  the {\em space of all $k$-jets of submanifolds of dimension}
$n$ of $W$. $J^k_n(W)$ has the following natural structures of differential fiber bundles:
$\pi_{k,s}:J^k_n(W)\to J^s_n(W)$, $s\le k$, with affine fibers
$J^k_n(W)_{\bar q}$, where $\bar q\equiv[N]^{k-1}_a\in J^{k-1}_n(W)$,
$a\equiv\pi_{k,0}(\bar q)$, with associated vector space $S^k(T^*_aN)\otimes\nu_a$,
$\nu_a\equiv T_aW/T_aN.$
For any $n$-dimensional submanifold $N\subset W$ one has the canonical
embedding $j^k:N\to J^k_n(W)$, given by $j^k:a\mapsto j^k(a)\equiv[N]^k_a$. We call
$j^k(N)\equiv N^{(k)}$ the {\em $k$-prolongation} of $N$.
In the following we shall also assume that there is a fiber
bundle structure on $W$, $\pi:W\to M$, where $\dim M=n.$ Then there exists a
canonical open bundle submanifold $J^k(W)$ of $J^k_n(W)$ that is called the
{\em $k$-jet space for sections} of $\pi$. $J^k(W)$ is diffeomorphic to the
$k$-jet-derivative space of sections of $\pi$, $J\mathcal{D}^k(W)$ \cite{PRA05}. Then,
for any section $s:M\to W$ one has the commutative diagram (\ref{commutative-diagram-derivative-space-jet-derivatives}), where $D^ks$ is the $k$-derivative of $s$ and $j^k(s)$ is the $k$-jet-derivative
of $s$. If $s(M)^{(k)}\subset J^k_n(W)$ is the $k$-prolongation of $s(M)\subset W$,
then one has $j^k(s)(M)\cong s(M)^{(k)}\cong s(M)\cong M$. Of course there are
also $n$-dimensional submanifolds $N\subset W$ that are not representable
as image of sections of $\pi$. As a consequence, in these cases, $N^{(k)}\cong N$ is not
representable in the form $j^k(s)(M)$ for some section $s$ of $\pi$.
The condition that $N$ is image of some (local) section $s$ of $\pi$ is equivalent
to the following local condition:
$s^*\eta\equiv s^*dx^1\wedge\dots\wedge dx^n\not=0$,
where $(x^\alpha,y^j)_{1\le\alpha\le n, 1\le j\le m}$, are fibered coordinates
on $W$, with $y^j$ vertical coordinates. In other words $N\subset W$ is locally representable by equations
$y^j=y^j(x^1,\dots,x^n)$. This is equivalent to saying that $N$ is transversal to the fibers of
$\pi$ or that the tangent space $TN$ identifies an horizontal
distribution with respect to the vertical one $vTW|_N$ of the
fiber bundle structure $\pi:W\to M$. Conversely, a completely integrable
$n$-dimensional horizontal distribution on $W$ determines a foliation of
$W$ by means of $n$-dimensional submanifolds that can be represented by
images of sections of $\pi$. The {\em Cartan distribution} of $J^k_n(W)$ is the distribution $\mathbf{E}^k_n(W)\subset TJ^k_n(W)$ generated by tangent spaces to the $k$-prolongation $N^{(k)}$ of $n$-dimensional submanifolds $N$ of $W$.

\begin{equation}\label{commutative-diagram-derivative-space-jet-derivatives}
  \xymatrix@C=2cm @R=2cm{J\mathcal{D}^k(W)\ar@{=}[r]^{\backsim}&J^k(W)\ar@{^{(}->}[r]&J^k_n(W)\ar@/^1pc/[dl]_{\pi_k}\\
&\ar@/^1pc/[ul]^{D^ks}M\ar[u]^{j^k(s)}&\\}
\end{equation}

\begin{theorem}[Metasymplectic structure of the Cartan distribution]\label{metasymplectic-structure}
There exists a canonical vector-fiber-valued $2$-form on the Cartan distribution $\mathbf{E}^k_n(W)$, called {\em metasymplectic structure} of $J^k_n(W)$.
\end{theorem}
\begin{proof}
The metasymplectic structure of the Cartan distribution $\mathbf{E}^k_n(W)\subset TJ^k_n(W)$ is a section $$\Omega_k:J^k_n(W)\to[S^{k-1}(\tau^*)\otimes\nu]\bigotimes\Lambda^2(\mathbf{E}^k_n(W)^*),$$
where $S^{k-1}(\tau^*)\equiv\bigcup_{q\in J^k_n(W)}S^{k-1}(\tau^*)_q$, with $S^{k-1}(\tau^*)_q\equiv S^{k-1}(T^*_aN)$, $\nu\equiv\bigcup_{q\in J^k_n(W)}\nu_q$, with $\nu_q\equiv(T_aW/T_aN)$, $[N]^k_a=q$, such that the following diagram
$$\xymatrix{S^{k-1}(\tau^*)_q\otimes\nu_q\ar@{=}[r]^{\backsim}\ar@{=}[d]^{\wr}& T_qJ^k_n(W)/\mathbf{E}^k_n(W)_q\ar@{=}[d]^{\wr}\\
T_{\bar q}J^{k-1}_n(W)/L_q\ar@{=}[r]^(0.35){\backsim}& \pi^{-1}_{k,k-1*}(T_{\bar q}J^{k-1}_n(W))/\pi^{-1}_{k,k-1*}(L_q)\\}$$
is commutative, for all $q\in J^k_n(W)$, $\bar q\equiv\pi_{k,k-1}(q)$, $a\equiv\pi_{k,0}(q)$, where $L_q\subset T_{\bar q}J^{k-1}_n(W)$ is the integral vector space canonically identified by $q$. Then, for the {\em metasymplectic structure} $\Omega_{\mathbf{E}^k_n(W)}$ of $\mathbf{E}^k_n(W)$ we have:
\begin{equation}\label{metasymplectic-structure-a}
\begin{array}{ll}
  \Omega_k(q)& \equiv\Omega_{\mathbf{E}^k_n(W)}(q)\in[T_qJ^k_n(W)/\mathbf{E}^k_n(W)_q]\bigotimes\Lambda^2(\mathbf{E}^k_n(W)^*_q) \\
  \\
  & \cong[S^{k-1}(\tau^*)_q\otimes\nu_q]\bigotimes\Lambda^2(\mathbf{E}^k_n(W)^*_q).\\
\end{array}
\end{equation}
More precisely $\Omega_k=d\omega_f|_{\mathbf{E}^k_n(W)}$, where $\omega_f=<\omega,f>=<f,(\phi^k)^*>\in\Omega^1(J^k_n(W))$ are the {\em Cartan forms} corresponding to smooth functions
$$f:J^k_n(W)\to\nu^k:=\bigcup_{q\in J^k_n(W)}\nu^k_q,\, \nu^k_q=T_{\bar q}J^{k-1}_n(W)/L_q.$$

$\phi^k$ is a canonical morphism of vector bundles over $J^k_n(W)$, defined by the exact sequence (\ref{exact-sequence-cartan-forms}).

\begin{equation}\label{exact-sequence-cartan-forms}
  \xymatrix{0\ar[r]&\mathbf{E}^k_n(W)\ar[rd]\ar[r]&TJ^k_n(W)\ar[d]\ar[r]^(0.6){\phi^k}&\nu^k\ar[ld]\ar[r]&0\\
  &&J^k_n(W)&&\\}
\end{equation}
For duality one has also the exact sequence (\ref{exact-sequence-cartan-forms-a}).

\begin{equation}\label{exact-sequence-cartan-forms-a}
  \xymatrix{0&\ar[l]\mathbf{E}^k_n(W)^*\ar[rd]&\ar[l]T^*J^k_n(W)\ar[d]&\ar[l]_(0.4){(\phi^k)^*}(\nu^k)^*\ar[ld]&\ar[l]0\\
  &&J^k_n(W)&&\\}
\end{equation}
Therefore we get also a smooth section $$\omega:J^k_n(W)\to\nu^k\bigotimes T^*J^k_n(W),$$ given by $<\omega,f>=<f,(\phi^k)^*>=f\circ \phi^k$, for any smooth section $f\in C^\infty((\nu^k)^*)$. It results
\begin{equation}\label{exact-sequence-cartan-forms-b}
  \mathbf{E}^k_n(W)=\bigcup_{f\in C^\infty((\nu^k)^*)}\ker(\omega_f).
\end{equation}

Furthermore, for any $\widetilde{q}\in \pi^{-1}_{k+1,k}(q)\subset J^{k+1}_n(W)$, $q=[N]^k_a\in  J^{k}_n(W)$, one has the following splitting:
\begin{equation}\label{exact-sequence-cartan-forms-c}
  \mathbf{E}^k_n(W)_q\cong L_{\widetilde{q}}\bigoplus [S^k(T^*_aN)\otimes\nu_a].
\end{equation}
The splitting (\ref{exact-sequence-cartan-forms-c}) allows us to give the following evaluation of $\Omega_k(q)(\lambda)$, for any $q\in J^k_n(W)$ and $\lambda\in S^{k-1}(T_aN)\otimes\nu_a^*$:
\begin{equation}\label{exact-sequence-cartan-forms-d}
\left\{
\begin{array}{ll}
  \Omega_k(q)(\lambda)(X,Y)=0,\,& \forall X,Y\in L_{\widetilde{q}},\, \pi_{k+1,k}(\widetilde{q})=q; \\
  \Omega_k(q)(\lambda)(\theta_1,\theta_2)=0,\, & \forall \theta_1,\theta_2\in S^k(T^*_aN)\otimes\nu_a; \\
  \Omega_k(q)(\lambda)(X,\theta)=<\lambda,X\rfloor\delta\theta>,\, & \forall X\in L_{\widetilde{q}},\, \theta\in S^k(T^*_aN)\otimes\nu_a,\\
  \end{array}
\right.
\end{equation}
where $\delta$ is the morphism in the exact sequence (\ref{delta-spencer-complex-a}).

If there is a fiber bundle structure $\alpha:W\to M$, $\dim M=n$, for the metasymplectic structure of $J\mathcal{D}^k(W)$ one has
$\Omega_k(q)\in \Lambda ^2(\mathbf{E}^k_n(W)^*_q)\bigotimes S^{k-1}(T^*_bM)\bigotimes
vT_aW$
with $a\equiv \pi _{k,0}(q)\in W$, $b\equiv \pi _k(q)\in M $.
If $\alpha $ is a trivial bundle
$\alpha :W\equiv M\times F\rightarrow M $, then
one has
$\Omega _k(q)\in \Lambda ^2(\mathbf{E}^k_n(W)^*_q)\bigotimes S^{k-1}(T^*_bM)
\bigotimes T_fF $, $\forall a\equiv (b,f) $.
\end{proof}

\begin{definition}
$\bullet$\hskip 2pt We say that vectors $X,Y\in \mathbf{E}^k_n(W)_q$ are in {\em involution} if  $$\Omega _k(q)(\lambda )(X,Y)=0,\, \forall \lambda \in S^{k-1}(T_aN)\bigotimes \nu ^*_a .$$

$\bullet$\hskip 2pt A subspace $P\subset \mathbf{E}^k_n(W)_q$ is called {\em isotropic} if any two vectors $X,Y\in P$ are in involution.

$\bullet$\hskip 2pt We say that a subspace $P\subset \mathbf{E}^k_n(W)_q$ is a
{\em maximal isotropic subspace} if $P$ is not a proper subspace of any other isotropic subspace.
\end{definition}
\begin{equation}\label{delta-spencer-complex-a}
\scalebox{0.9}{$\xymatrix{0\ar[d]\\[S^m(T^*_aN)\bigotimes\nu_a]\ar[d]^(0.35){\delta}\\
T^*_aN\bigotimes [S^{m-1}(T^*_aN)\otimes\nu_a]\ar[d]^(0.45){\delta}\\
\Lambda^2(T^*_aN)\bigotimes[S^{m-2}(T^*_aN)\otimes\nu_a]\ar[d]^(0.5){\delta}\\
\cdots\ar[d]^(0.3){\delta}\\
\Lambda^n(T^*_aN)\bigotimes [S^{m-n}(T^*_aN)\otimes\nu_a]\ar[d]\\
0\\}$}
\end{equation}
\begin{theorem}[Structure of maximal isotropic subspaces]\label{structure-maximal-isotropic-subspaces}
Any maximal isotropic subspace $P\subset \mathbf{E}^k_n(W)_q$ is one tangent at $q=[N]^k_a$ to a maximal integral manifold $V$ of $\mathbf{E}^k_n(W)$. These are of dimension $m\binom{p+k-1}{k}+n-p$, such that $n-p=\dim(\pi_{k,0}{}_*(T_qV))\le \dim T_aN=n$. Then one says that $V$ is of {\em type} $n-p$. In particular if $p=0$, then $L_{\widetilde{q}}\cong T_qV\cong T_aN$. In the {\em exceptional case}, i.e., $m=n=1$, maximal integral manifolds are of dimension $1$ having eventual subsets belonging to the fibers of $\pi_{k,k-1}:J^k_n(W)\to J^{k-1}_n(W)$.
\end{theorem}
\begin{proof}
The degeneration subspace of $\Omega_k(q)(\lambda)$, for any $\lambda\in S^{k-1}(T_aN)\otimes\nu_a^*$, is the subaspace $P\subset \mathbf{E}^k_n(W)_q$ given in (\ref{structure-maximal-isotropic-subspaces}).
\begin{equation}\label{structure-maximal-isotropic-subspaces}
P\equiv \left\{ <x+ \theta >\vert   x\in Ann(\Xi)\subset T_aN ,\,
\theta \in S^k(\Xi)\bigotimes \nu _a\subset S^k(T^*_aN)\bigotimes \nu _a\right\},
\end{equation}
where $\Xi$ is a $p$-dimensional subspace of $T^*_aN $.

Let, now, $N\subset W$ be a $n$-dimensional
submanifold of $W$ and let $N_0\subset N$ be
a submanifold in $N$. Set
$$N^{(k)}_0(N)\equiv \left\{  q\in J^k_n(W)\, \vert\,  \pi _{k,k-1}(q)\in N^{(k-1)}_0 ,\, L_q\supset T_{\pi _{k,k-1}(q)}N^{(k-1)}_0\right\} $$
where $N^{(k-1)}_0\equiv \left\{ [N]^{k-1}_a\, \vert\,   a\in N_0\right\}\subset J^{k-1}_n(W) $.
Then the tangent planes to $N^{(k)}_0(N)$ coincide with the maximal involutive subspaces described in (\ref{structure-maximal-isotropic-subspaces}). Therefore,
$N^{(k)}_0(N)$ is a maximal integral manifold of the Cartan distribution.
\end{proof}
\begin{equation}\label{delta-spencer-complex}
\scalebox{0.9}{$\xymatrix@R=1cm{0\ar[d]\\
 g_m(q)\ar[d]^(0.35){\delta}\\
 T^*_aN\bigotimes g_{m-1}(q)\ar[d]^(0.45){\delta}\\
\Lambda^2(T^*_aN)\bigotimes g_{m-2}(q)\ar[d]^(0.5){\delta}\\
\ar[d]^(0.3){\delta}\cdots\\
\ar[d]\Lambda^n(T^*_aN)\bigotimes g_{m-n}(q)\\
0\\}$}
\end{equation}
\begin{definition}[Partial differential equation for submanifolds]\label{definition-pde}
A {\em partial differential equation} (PDE) for $n$-dimensional submanifolds
of $W$ is a submanifold $E_k\subset J^k_n(W)$.\footnote{In this paper, for sake of simplicity, we shall consider only smooth PDEs. For information on the geometry of singular PDEs, see the following references \cite{AG-PRA2,PRA3,PRA9,PRA11,PRA14}.} A ({\em regular}) {\em solution}
of $E_k$ is a (regular) solution of $J^k_n(W)$ that is contained into
$E_k$. In particular, if $E_k\subset J^k(W)\subset J^k_n(W)$ we can talk
about PDE for sections of $\pi:W\to M$.
The {\em prolongation of order $l$} of $E_k\subset J^k_n(W)$  is the subset
$(E_k)_{+l}\subset J^{k+l}_n(W)$ defined by
$(E_k)_{+l}\equiv J^l_n(E_k)\cap J^{k+l}_n(W)$.
A PDE $E_k\subset J^k_n(W)$ is called {\em formally integrable} if for all
$l\ge 0$ the prolongations $(E_k)_{+l}$ are smooth submanifolds and the projections
$\pi_{k+l+1,k+l}:(E_k)_{+(l+1)}\to (E_k)_{+l}$, $\pi_{k,0}:E_k\to W$ are smooth
bundles. The {\em symbol} of the PDE $E_k\subset J^k_n(W)$ at the point
$q\equiv[N]^k_a\in E_k$ is defined to be the following subspaces:
$g_k(q)\equiv T_{q}(E_k)\cap T_{q}(F_{\bar q})$, where $\bar q\equiv\pi_{k,k-1}(q)$, and $\pi^{-1}_{k,k-1}(\bar q)=F_{\bar q}\subset J^k_n(W)$.
Using the affine structure on the fibre $F_{\bar q}$, we can identify the
symbol $g_k(q)$ with a subspace in $S^k(T^*_aN)\otimes\nu_a$:
$g_k(q)\subset S^k(T^*_aN)\otimes\nu_a$. Suppose that all prolongations
$(E_k)_{+l}$ are smooth manifolds, then their symbols at points
$\breve q\equiv[N]^{k+l}_a$ are $l$th prolongations of the symbol $g_k(q)$, hence
$g_{k+l}(\breve q)=g_{k+l}(q)\subset S^{k+l}(T^*_aN)\otimes \nu_a$
and
$\delta(g_{k+l}(\breve q))\subset g_{k+(l-1)}(q)\otimes T^*_aN$, $l=1,2,\dots$
where by
$\delta:S^{k+l}(T^*_aN)\otimes \nu_a\to T^*_aN\otimes
S^{k+l-1}(T^*_aN)\otimes \nu_a$
we denote {\em $\delta$-Spencer operator}. Therefore, at each point $q\in E_k$
the $\delta$-Spencer complex is defined, where $m\ge k$. We denote by $H^{m-j,j}(E_k,q)$ the cohomologies
of this complex at the term  $\Lambda^j(T^*_aN)\otimes g_{m-j}(q)$. They are
called {\em $\delta$-Spencer cohomologies} of PDE at the point $q\in E_k$.
We say that $g_k$ is {\em involutive} if the sequences {\em(\ref{delta-spencer-complex})}
are exact and that $g_k$ is {\em $r$-acyclic} if $H^{m-j,j}(E_k,q)=0$ for
$m-j\ge k$, $0\le j\le r$.
If $E_k\subset J^k_n(W)$ is a {\em $2$-acyclic} PDE , i.e.,
$H^{j,i}(E_k,q)=0$, $\forall q\in E_k$, $ 0\le j\le 2$, $m-j\ge k$,
and $\pi_{k+1,k}:E_k^{(1)}\to E_k$, $\pi_{k,0}:E_k\to W$ are smooth bundles, then
$E_k$ is {\em formally integrable}.
\end{definition}
\begin{definition}
We say that $E_k\subset J^k_n(W)$ is {\em completely integrable} if for any point $q\in E_k$, passes a (local) solution of $E_k$, hence a $n$-dimensional manifold $V\subset E_k$, with $q\in V$ and $V=N^{(k)}$. This implies that the following sequence
$$\xymatrix@C=2cm{(E_k)_{+r}\ar[r]^{\pi_{k+r,k+r-1}}&E_{k+r-1}\ar[r]&0}$$
is exact for any $r\ge 1$. (This is equivalent to say that $\pi_{k+r,k+r-1}|_{(E_k)_{+r}}$ is surjective.
\end{definition}
\begin{proposition}\label{analytical-category}
In the category of analytic manifolds, (i.e., manifolds of class $C^\omega$), the formal integrability implies the complete integrability.
\end{proposition}

\begin{definition}
A {\em Cartan connection} on $E_k$ is a $n$-dimensional subdistribution
$\mathbf{H}\subset \mathbf{E}_k$ such that $T(\pi_{k,k-1})(\mathbf{H}_q)=L_q\equiv T_{\pi_{k,k-1}(q)}N^{(k-1)}$, $[N]^k_q\equiv q,\quad\forall q\in E_k$.\footnote{As $\dim(L_q)=n=\dim\mathbf{H}_q$ then there exists a $n$-dimensional
submanifold $X\subset W$ such that $T_qX^{(k)}=\mathbf{H}_q$,
with $[X]^k_a=q$, $[X]^{k-1}_a=[N]^{k-1}_a$, $T_{\pi_{k,k-1}(q)}X^{(k-1)}=L_q$.}
 We call {\em curvature} of the Cartan connection
$\mathbf{H}$ on $E_k\subset J^k_n(W)$ the field of geometric objects on $E_k$:
\begin{equation}\label{curvature-cartan-connection}
\begin{array}{ll}
  \Omega_{\mathbf{H}}: & q\mapsto\Lambda^2(\mathbf{H}^*_q)\bigotimes[S^{k-1}(T_aN)\otimes\nu_a^*/Ann(g_{k-1})]^*\\
  & \cong\Lambda^2(T^*_aN)\bigotimes[S^{k-1}(T_aN)\otimes\nu_a^*/Ann(g_{k-1})]^*\\
\end{array}
\end{equation}
obtained by restriction on $\mathbf{H}$ of the metasymplectic structure on the
distribution $\mathbf{E}^k_n$.
\end{definition}

\begin{proposition}
In any {\em flat Cartan connection} $\mathbf{H}\subset\mathbf{E}_k$, i.e., a Cartan connection having zero curvature: $\Omega_{\mathbf{H}}=0$, any two vector $X,Y\in \mathbf{H}_q$, $q\in E_k$ are in involution.
\end{proposition}

\begin{definition}
Let us assume that $(E_k)_{+1}\to E_k$ is a smooth subbundle of
$J^{k+1}_n(W)\to J^k_n(W)$. Then any section $\rceil:E_k\to(E_k)_{+1}$ is called a {\em Bott connection}.
\end{definition}

\begin{theorem}

{\rm 1)} A Cartan connection $\mathbf{H}$ is a Bott connection iff
$\Omega_{\mathbf{H}}=0$.\footnote{If $(E_k)_{+1}\to E_k$ is a smooth subbundle of
$J^{k+1}_n(W)\to J^k_n(W)$ then a flat Cartan connection is also an involutive distribution. On the other hand a Bott connection identifies an involutive distribution iff it is a flat connection. (For more details on $(k+1)$-connections on $W$, see \cite{PRA1}.)}

{\rm 2)} A Cartan connection $\mathbf{H}$ gives a splitting of the Cartan distribution
$$\mathbf{E}^k_n\cong g_k\bigoplus\mathbf{H}.$$
Two Cartan connections $\mathbf{H}$, $\mathbf{H}'$ on $E_k$ identify a field of
geometric objects $\lambda$ on $E_k$ called {\em soldering form}: $\lambda\equiv\lambda_{\mathbf{H},\mathbf{H}'}:E_k\to\mathbf{H}^*\bigotimes g_k$, $\lambda(q)\in T_a^*N\otimes g_k(q)$. One has:

$\bullet$\hskip 2pt $\Omega_{\mathbf{H}'}=\Omega_{\mathbf{H}}+\delta\lambda$.

$\bullet$\hskip 2pt {\em(Bianchi identity)} $\delta\Omega_{\mathbf{H}}=0$,
$$\Omega_{\mathbf{H}}(q)\hskip 2pt \hbox{\rm mod}\hskip 2pt
         \delta(T_a^*N\otimes g_k(q))\in H^{k-1,2}(E_k)_q.$$ We call such $\delta$-cohomology class of $\Omega_{\mathbf{H}}$ the {\em Weyl tensor}
of $E_k$ at $q\in E_k$: $W_k(q)\equiv[\Omega_{\mathbf{H}}(q)]$.
Then, there exists a point $u\in(E_k)_{+1}$ over $q\in E_k$ iff $W_k(q)=0$.

{\em 3)} Suppose that $g_{k+1}$ is a vector bundle over $E_k\subset J^k_n(W)$.
Then if the Weyl tensor $W_k$ vanishes the projection $\pi_{k+1,k}:(E_k)_{+1}\to E_k$
is a smooth affine bundle.

{\em 4)} If $g_{k+l}$ are vector bundles over $E_k$ and $W_{k+l}=0$, $l\ge 0$,
then $E_k$ is formally integrable.

{\em 5)} If the system $E_k$ is of finite type, i.e., $g_{k+l}(q)=0$,
$\forall q\in E_k$, $l\ge l_0$, then $W_{k+l}=0$, $0\le l\le l_0$, is a
sufficient condition for integrability.
\end{theorem}

\begin{theorem}
Given a Cartan connection $\mathbf{H}$ on $E_k$, for any regular solution
$N^{(k)}\subset E_k$ we identify a section
${}_{\mathbf{H}}\nabla\in C^\infty(T^*N\bigotimes g_k)$
called {\em covariant differential} of $\mathbf{H}$ of the solution $N$. Furthermore,
for any vector field $\zeta:N\to TN$ we get a section
${}_{\mathbf{H}}\nabla\zeta\in C^\infty(g_k|_{N^{(k)}}).$
\end{theorem}

\begin{theorem}[Characteristic distribution of PDE]
Let $E_k\subset J^k_n(W)$ be a PDE such that $(E_k)_{+1}\to E_k$ is a smooth subbundle of $J^{k+1}_n(W)\to J^k_n(W)$. Then for any $\widetilde{q}\in(E_k)_{+1}$ the set $\mathbf{Char}(E_k)_q$ of vectors in the splitting $(\mathbf{E}_k)_q\cong L_{\widetilde{q}}\bigoplus(g_k)_q$, $\zeta=v+\theta$, such that $v\rfloor\delta(\theta)=0$, for any $\theta\in(g_k)_q$ is called the {\em space of characteristic vectors} at $q\in E_k$.
$\mathbf{Char}(E_k)$ is an involutive substribution of the Cartan distribution $\mathbf{E}_k$.

$\bullet$\hskip 2pt $\mathbf{Char}(E_k)=\mathbf{E}_k\bigcap \mathfrak{s}(E_k)$, where $\mathfrak{s}(E_k)$is the space of infinitesimal simmetries of $E_k$, namely the set of vector field on $E_k$  whose flows preserve the Cartan distribution.
\end{theorem}
\begin{proof}
See \cite{PRA1}.
\end{proof}

\begin{definition}
We call a PDE $E_k\subset J^k_n(W)$ {\em degenerate} at the point $q\in E_k$ if there is a $p$-dimensional ($0<p\le n$), subspace $\Xi_q\subset T_a^*N$, such that
$$(g_k)_q\subset [S^k(\Xi_q)\bigotimes\nu_a].$$
\end{definition}

\begin{theorem}
$\mathbf{Char}(E_k)_q\not= 0$ iff $E_k$ is a degenerate PDE at the point $q\in E_k$. The subspace
$$\Xi_q=Ann((\pi_{k,0})_*(\mathbf{Char}(E_k)_q))$$
is the subspace of degeneration of $E_k$ at the point $q\in E_k$.

$\bullet$\hskip 2pt Let $E_k\subset J^k_n(W)$ be a PDE such that the following conditions hold:

{\em (i)} $\pi_{k+1,k}:(E_k)_{+1}\to E_k$ and $\pi_{k,k-1}:E_k\to J^{k-1}_n(W)$ are smooth bundles:

 {\em (ii)} $\Xi=\bigcup_{q\in E_k}\Xi_q$ is a smooth vector bundle, where $\Xi_q$ is a space of degeneration of $E_k$ at the point $q\in E_k$. Then, $\mathbf{Char}(E_k)$ is a smooth distribution on $E_k$. and solutions of $E_k$ can be formulated by the method of characteristics.\footnote{In other words the method of characteristics allows us to solve Cauchy problems in $E_k$, namely to build a solution $V$ containing a fixed $(n-1)$-dimensional integral manifold $N_0$: $N_0\subset V$. In fact if $\zeta:E_k\to TE_k$ is a characteristic vector field of $E_k$, transverse to $N_0$, then $V=\bigcup_t\phi_t(N_0)$ is a solution of $E_k$, if $\partial \phi=\zeta$.}
\end{theorem}

In this section we shall classify global singular solutions of PDEs by means of suitable bordism groups.

\begin{definition}[Generalized singular solutions of PDE]\label{generalized-singular-solutions-of-pde}
Let $E_k\subset J^k_n(W)$ be a PDE.
We call {\em bar singular chain complex}, {\em with coefficients
into an abelian group $G$}, of $E_k$ the chain complex:
 $$\{\bar C_p(E_k;G),\bar\partial\},$$
 where $\bar C_p(E_k;G)$ is the
$G$-module of formal linear combinations, with coefficients in
$G$, $\sum \lambda_i c_i$, where $c_i$ is a singular $p$-chain
$f:\bigtriangleup ^p\to E_k$ that extends on a neighborhood
$U\subset\mathbb{R}^{p+1}$, such that $f$ on $U$ is differentiable
and $Tf(\bigtriangleup^p)\subset\mathbf{E}_k$. Denote by $\bar
H_p(E_k;G)$ the corresponding homology ({\em bar singular homology
with coefficients in $G$}) of $E_k$.

A {\em $G$-singular $p$-dimensional integral manifold} of $E_k\subset
J^k_n(W)$, is a bar singular $p$-chain $V$ with $p\le n$, and coefficients into an abelian group $G$, such that $V\subset
E_k$.

Set $\bar B_\bullet(E_k;G)\equiv\IM(\bar\partial)$, $\bar Z_\bullet(E_k;G)\equiv
\ker(\bar\partial)$. Therefore, one has the exact
commutative diagram {\em(\ref{exact-commutative-diagram-bordism-g-singular-groups})}.
\begin{equation}\label{exact-commutative-diagram-bordism-g-singular-groups}
\xymatrix{&&0\ar[d]&0\ar[d]&&\\
               &0\ar[r]&\bar B_\bullet(E_k;G)\ar[r]\ar[d]&\bar Z_\bullet(E_k;G)\ar[r]\ar[d]&\bar H_\bullet(E_k;G)\ar[r]&0\\
&&\bar C_\bullet(E_k;G)\ar[d]\ar@{=}[r]& \bar C_\bullet(E_k;G)\ar[d]&&\\
0\ar[r]& {}^G\Omega^{E_k}_{\bullet,s}\ar[r]&\bar Bor_\bullet(E_k;G)\ar[d]\ar[r]&\bar Cyc_\bullet(E_k;G)\ar[d]\ar[r]&0&\\
&&0&0&&\\}
\end{equation}
In Tab. \ref{legenda-for-commutative-exact-diagram-singular-g-solution} are given some more explicit properties about the symbols involved in {\em(\ref{exact-commutative-diagram-bordism-g-singular-groups})}.
\end{definition}

\begin{table}[t]
\caption{Legenda for the commutative exact diagram (\ref{exact-commutative-diagram-bordism-g-singular-groups}).}
\label{legenda-for-commutative-exact-diagram-singular-g-solution}
\centerline{\scalebox{0.9}{$\begin{tabular}{|l|l|l|}
\hline
\hfil{\rm{\footnotesize Name}}\hfil&\hfil{\rm{\footnotesize Definition}}\hfill&\hfil{\rm{\footnotesize Properties}}\hfill\\
\hline\hline
\hfil{\rm{\footnotesize Bordism group}}\hfil&\hfil{\rm{\footnotesize $\bar Bor_\bullet(E_k;G)$}}\hfill&\hfil{\rm{\footnotesize $b\in{}^G[a]_{E_k}\in\bar Bor_\bullet(E_k;G)\Rightarrow\exists c\in \bar C_\bullet(E_k;G): \bar\partial c=a-b$}}\hfill\\
\hline
\hfil{\rm{\footnotesize Cyclism group}}\hfil&\hfil{\rm{\footnotesize $\bar Cyc_\bullet(E_k;G)$}}\hfill&\hfil{\rm{\footnotesize $b\in{}^G[a]_{E_k}\in\bar Cyc_\bullet(E_k;G)\Rightarrow \bar\partial(a-b)=0$}}\hfill\\
\hline
\hfil{\rm{\footnotesize Closed bordism group}}\hfil&\hfil{\rm{\footnotesize ${}^G\Omega^{E_k}_{\bullet,s}$}}\hfill&\hfil{\rm{\footnotesize $b\in{}^G[a]_{E_k}\in{}^G\Omega^{E_k}_{\bullet,s}\Rightarrow
\left\{
\begin{array}{l}
 \bar\partial a=\bar\partial b=0\\
a-b=\bar\partial c\\
\end{array}\right\}$}}\hfill\\
\hline
\end{tabular}$}}
\end{table}

\begin{theorem}[Integral singular bordism groups of PDE]\label{integral-singular-bordism-groups-of-pde}

$\bullet$\hskip 2pt One has the following canonical isomorphism:
$${}^G\Omega^{E_k}_{\bullet,s}\cong\bar H_\bullet(E_k;G).$$
$\bullet$\hskip 2pt If ${}^G\Omega^{E_k}_{\bullet,s}=0$ one has: $\bar Bor_\bullet(E_k;G)\cong\bar Cyc_\bullet(E_k;G)$.

$\bullet$\hskip 2pt If $\bar Cyc_\bullet(E_k;G)$ is a free $G$-module, then the bottom horizontal exact sequence, in above diagram, splits and one has the isomorphism:
$$\bar Bor_\bullet(E_k;G)\cong {}^G\Omega( E_k)_{\bullet,s}\bigoplus\bar Cyc_\bullet(E_k;G).$$
\end{theorem}
\begin{remark}
By considering the dual complex
\begin{equation}\label{dual-complex}
  \{\bar C^p(E_k;G)\equiv
Hom_{\mathbb{Z}}(\bar C_p(E_k;{\mathbb{Z}});G), \bar\delta\}
\end{equation}

and $\bar H^p(E_k;G)$, the associated
homology spaces ({\em bar singular cohomology}, {\em with
coefficients into $G$} of $E_k$), we can talk also of singular co-bordism groups with coefficients in $G$. These are important objects, but in this paper we will skip on these aspects.
\end{remark}
\begin{definition}
A {\em $G$-singular $p$-dimensional quantum manifold} of $E_k$ is a bar
singular $p$-chain $V\subset J^k_n(W)$, with $p\le n$, and coefficients into an abelian group
$G$, such that $\partial V\subset E_k$. Let us denote by
${}^G\Omega_{p,s}(E_k)$ the corresponding (closed) bordism
groups in the singular case. Let us denote also by
${}^G[N]_{\overline{E_k}}$ the equivalence classes of quantum
singular bordisms respectively.\footnote{These bordism groups can be called also $G$-singular $p$-dimensional integral bordism groups {\em relative} to $E_k\subset J^k_n(W)$. They play an important role in PDE algebraic topology. For more details see Refs. \cite{PRA01,PRA1,PRA7,PRA8,PRA9,PRA10,PRA11}.}
\end{definition}

\begin{remark}
Let us emphasize that a $G$-singular solution $V\subset E_k$ can be written as a $n$-chain $V=\sum _ia^i u_i$, where
$a_i\in G$ and $u_i:\Delta^n\to E_k$, such that $u_i(\Delta^n)$ is an integral manifold of $E_k$.\footnote{In such a category can be considered also so-called {\em neck-pinching singular solutions} that are very important whether from a theoretical point of view as well in applications. (See, e.g., \cite{PRA15,PRA16}.)} In particular a $G$-singular solution $V$ of $E_k$ can have tangent spaces $T_qV$ is some points $q\in V$ such that $T_qV$ is a {\em $n$-dimensional integral plane}, i.e., an $n$-dimensional subspace of $(\mathbf{E}_k)_q\subset T_q E_k$ of the type $L_{\widetilde{q}}$, for some $\widetilde{q}\in (E_k)_{+1}$, or admitting the splitting
$$T_qV=V_q^k\bigoplus V_q^0$$
where $V_q^k=T_qV\bigcap(g_k)_q\subset V_q\bigcap[S^k(T^*_aN)\bigotimes\nu_a]$ and $V_q^0\subset L\widetilde{q}$, $V_q^0\cong(\pi_{k,0}(V_q))\subset T_aN$, $\dim V_q^0={\rm type}(V)=n-p$. $(g_k)_q$ is the unique maximal isotropic subspace of dimension equal to $m\binom{p+k-1}{k}$ (and type $0$). Therefore, under the condition (\ref{condition-singularity}).

\begin{equation}\label{condition-singularity}
  m\binom{p+k-1}{k}\ge n
\end{equation}
a singular solution of $E_k$ can contain  pieces of type $0$. We say that a singular solution is {\em completely degenerate} if it is an integral $n$-chain of type $0$, namely completely contained in the symbol $(g_k)_q$, for some $q\in E_k$. In general a singular solution can contain completely degenerate pieces.  When the set $\Sigma(V)\subset V$ of singular points of a singular solution $V\subset E_k$, is nowhere dense in $V$, therefore $\dim \Sigma(V)<n$, then we say that in $V$ there are {\em Thom-Boardman singularities}. In such points $q\in V$ one has $\dim[T_qV\bigcap(g_k)_q]=p$, with $0<p<n$. This is equivalent to state that $\dim[(\pi_{k,0})_*(T_qV)]=n-p$, or that $q$ is a point of {\em Thom-Boardman-degeneration}.
Finally when $\Sigma(V)=\varnothing$, and there are not completely degenerate points in $V$, we say that $V$ is a {\em regular solution}. In such a case $V$ is diffeomorphic to its projection $X=\pi_{k,0}(V)\subset W$, or equivalently $\pi_{k,0}|_V:V\to W$ is an embedding.
\end{remark}

\begin{theorem}[Cauchy problems in PDE]\label{singular-cauchy-problems-in-pde}
If $E_k$ is a completely integrable PDE, and $\dim (g_k)_{+1}\ge n$, given a $(n-1)$-dimensional regular integral manifold $N$, contained in $E_k$, there exists a solution $V\subset E_k$, such that $V\supset N$.
\end{theorem}
\begin{proof}
In fact, since $N$ is regular, it identifies a $(n-1)$-manifold in $W$, say $N_0\subset W$. Let  $Y\subset W$ be a $n$-dimensional manifold containing $N_0$. Then taking into account that $E_k$ is completely integrable, we can assume that the $(k+1)$-prolongation $Y^{(k+1)}\subset J^{k+1}_nW$ of $Y$ is such that $Y^{(k+1)}\bigcap(E_k)_{+1}=N_0^{(1)}$, namely it coincides with an $(n-1)$-dimensional integral manifold that projects on $E_k$. We call $N_0^{(1)}$ the first prolongation of $N_0$. Now taking into account that $(E_k)_{+1}$ is the strong retract of $J^{k+1}_n(W)$, we can retract map $Y^{(k+1)}$ into $(E_k)_{+1}$, via the retraction, obtaining a solution $V'\subset(E_k)_{+1}$ of $(E_k)_{+1}$ passing for $N_0^{(1)}$. By projecting $V'$ into $E_k$, we obtain a solution $V$ containing $N$. Since $\dim (g_k)_{+1}\ge n$, the solution $V'$ does not necessitate to be regular, but can have singular points.
\end{proof}
\begin{example}
Let $E_2\subset J\mathcal{D}^2(W)$, be an analytic dynamic equation of a rigid system with $n$-degree of freedoms. Let $\{t,q^i.\dot q^i,\dot q^i\}$ be local coordinates on $J\mathcal{D}^2(W)$. Such an equation is completely integrable. A Cauchy problem there is encoded by a point $q_0\in E_2$, hence for that point pass an unique solution V, i.e., an integral curve contained into $E_2$. Let us, however, try to apply the proceeding of the proof of Theorem \ref{singular-cauchy-problems-in-pde}. This is strictly impossible ! In fact the symbol of such an equation is necessarily zero: $\dim(g_2)_q=0$, for any $q\in E_2$.\footnote{In general such dynamical equations have zero symbol since they are encoded by $n$ analytic differential equations of the second order, where $n$ is the degree of freedoms.} On the other hand we can consider a point $\widetilde{q}_0$ belonging to $(E_2)_{+1}$ and such that $\pi_{3,2}(\widetilde{q}_0)=q_0$, and $\pi_{2,0}(\widetilde{q}_0)=a\in W$, and we can assume that there exists an integral curve $Y\subset J\mathcal{D}^3(W)$ passing for $\widetilde{q}_0$, but when we retract such a curve into $(E_2)_{+1}$, we get the unique curve $\overline{\Gamma}$ passing for $\widetilde{q}_0$ contained into $(E_2)_{+1}$. This curve does not necessarily pass for the point $\bar q_0=V^{(1)}\bigcap\pi^{-1}_{3,2}(q_0)$, since the first prolongation $V^{(1)}$ of $V$ does not necessarily coincide with $\overline{\Gamma}$. Thus the proceeding considered in the proof of Theorem \ref{singular-cauchy-problems-in-pde} does not apply to PDEs (or ODEs), having zero symbols $g_2=0$. In other words, for such PDEs, despite $\pi_{2,0}(q_0)=\pi_{2,0}(\overline{q}_0)=a\in W$,
 we cannot connect two regular solutions corresponding to two different initial conditions $q_0$ and $\overline{q}_0$, with a completely degenerate piece, or a Thom-Boardman-singular piece. However, a more general concept of solutions can be considered also when $g_k=0$. In fact {\em weak solutions} allow include solutions with discontinuity points.\footnote{It is worth to emphasize that weak solutions can be considered equivalent to solutions having completely degenerated pieces, in fact their projections on the configuration space $W$ are the same. However weak solution can exist also with trivial symbol $g_k=0$, instead solutions with completely degenerated pieces can exist only if $\dim g_k\ge n$. Furthermore, under this circumstance, namely under condition (\ref{condition-singularity}), a {\em continuous weak solution}, i.e., a weak solution having completely degenerate pieces, can be deformed into solutions with Thom-Boardman singular points.}
\end{example}
\begin{remark}
Weak solutions are of great importance and must be included in a geometric theory of PDE's too.
\end{remark}

\begin{definition}
Let $\Omega_{n-1}^{E_k}$, (resp. $\Omega_{n-1,s}^{E_k}$, resp.
$\Omega_{n-1,w}^{E_k}$), be the integral bordism group for
$(n-1)$-dimensional smooth admissible regular integral manifolds
contained in $E_k$, bounding smooth regular integral
manifold-solutions,\footnote{This means that
$N_1\in[N_2]\in\Omega_{n-1}^{E_k}$, iff
$N_1^{(\infty)}\in[N^{(\infty)}_2]\in\Omega_{n-1}^{E_{\infty}}$.
(See Refs.\cite{PRA4, PRA14} for notations.)} (resp. piecewise-smooth or singular
solutions, resp. singular-weak solutions), of $E_k$.
\end{definition}

\begin{theorem}\label{exact-commutative-diagram-relating-smooth-seingular-weak-integral-bordism-groups}
Let $\pi:W\to M$ be a fiber bundle with $W$ and $M$ smooth
manifolds, respectively of dimension $m+n$ and $n$. Let $E_k\subset
J^k_n(W)$ be a PDE for $n$-dimensional submanifolds of $W$.  One has
the following exact commutative diagram relating the groups
$\Omega_{n-1}^{E_k}$, $\Omega_{n-1,s}^{E_k}$ and
$\Omega_{n-1,w}^{E_k}$:

\begin{equation}\label{commutative-diagram-relating-bordism-groups-in-pde}
\xymatrix{
&0\ar[d]&0\ar[d]&0\ar[d]&\\
0\ar[r]&K^{E_k}_{n-1,w/(s,w)}\ar[d]\ar[r]&
K^{E_k}_{n-1,w}\ar[d]\ar[r]&K^{E_k}_{n-1,s,w}\ar[d]\ar[r]&0\\
0\ar[r]&K^{E_k}_{n-1,s}\ar[d]\ar[r]&
\Omega^{E_k}_{n-1}\ar[d]\ar[r]&\Omega^{E_k}_{n-1,s}\ar[d]\ar[r]&0\\
&0\ar[r]&\Omega^{E_k}_{n-1,w}\ar[d]\ar[r]&\Omega^{E_k}_{n-1,w}\ar[d]\ar[r]&0\\
&&0&0&}
\end{equation}
and the canonical isomorphisms reported in {\em(\ref{canonical-isomorphisms})}.
\begin{equation}\label{canonical-isomorphisms}
\left\{ \begin{array}{l}
 K^{E_k}_{n-1,w/(s,w)}\cong K^{E_k}_{n-1,s} \\
 \Omega^{E_k}_{n-1}/K^{E_k}_{n-1,s}\cong
\Omega^{E_k}_{n-1,s}\\
\Omega^{E_k}_{n-1,s}/K^{E_k}_{n-1,s,w}\cong\Omega^{E_k}_{n-1,w}\\
\Omega^{E_k}_{n-1}/K^{E_k}_{n-1,w}\cong\Omega^{E_k}_{n-1,w}.\\
 \end{array}\right.
\end{equation}
$\bullet$\hskip 2pt In particular, for $k=\infty$, one has the canonical
isomorphisms reported in {\em(\ref{canonical-isomorphisms-a})}.
\begin{equation}\label{canonical-isomorphisms-a}
\left\{ \begin{array}{l}
K^{E_\infty}_{n-1,w}\cong K^{E_\infty}_{n-1,s,w}\\
K^{E_\infty}_{n-1,w/(s,w)}\cong K^{E_\infty}_{n-1,s}\cong 0\\
\Omega^{E_\infty}_{n-1}\cong \Omega^{E_\infty}_{n-1,s}\\
\Omega^{E_\infty}_{n-1}/K^{E_\infty}_{n-1,w}\cong\Omega^{E_\infty}_{n-1,s}/K^{E_\infty}_{n-1,s,w}\cong\Omega^{E_\infty}_{n-1,w}.\\
\end{array}\right.
\end{equation}

$\bullet$\hskip 2pt If $E_k$ is formally integrable then
one has the isomorphisms reported in {\em(\ref{canonical-isomorphisms-b})}.
\begin{equation}\label{canonical-isomorphisms-b}
\Omega^{E_k}_{n-1}\cong\Omega^{E_\infty}_{n-1}\cong\Omega^{E_\infty}_{n-1,s}.
\end{equation}
\end{theorem}

\begin{proof}
The proof follows directly from the definitions and standard
results of algebra. (For more details see Refs. \cite{PRA1,PRA6},)
\end{proof}

\begin{theorem}\label{formal-integrability-integral-bordism-groups}
Let us assume that $E_k$ is formally integrable and completely
integrable, and such that $\dim E_k\ge 2n+1$. Then, one has the canonical isomorphisms reported in {\em(\ref{canonical-isomorphisms-c})}.
\begin{equation}\label{canonical-isomorphisms-c}
\Omega^{E_k}_{n-1,w}\cong\bigoplus_{r+s=n-1}H_r(W;{\mathbb Z
}_2)\otimes_{{\mathbb Z}_2}\Omega_s
\cong\Omega^{E_k}_{n-1}/K^{E_k}_{n-1,w}\cong
\Omega^{E_k}_{n-1,s}/K^{E_k}_{n-1,s,w}.
\end{equation}
where $\Omega_s$ denotes the $s$-dimensional un-oriented smooth bordism group.

$\bullet$\hskip 2pt Furthermore, if $E_k\subset J^k_n(W)$, has non zero symbols: $g_{k+s}\not=0$,
$s\ge 0$, (this excludes that can be $k=\infty$), then $K^{E_k}_{n-1,s,w}=0$, hence $\Omega^{E_k}_{n-1,s}\cong\Omega^{E_k}_{n-1,w}$.
\end{theorem}

\begin{proof}
It follows from above theorem and results in \cite{PRA4}.
Furthermore, if $g_{k+s}\not=0$, $s\ge 0$, we can always connect two
branches of a weak solution with a singular solution of $E_k$. (For more details see \cite{PRA4}.)
\end{proof}

\section{\bf Maslov index in PDEs and Lagrangian bordism groups}\label{sec:4}

In order to consider ``Maslov index" canonicallay associated to PDEs, we follow a strategy to recast Arnold-Kashiwara-Thomas algebraic approach, resumed in section 2, by substituting the Grassmannian of Lagrangian subsapces with the Grassmannian of $n$-dimensional integral planes, namely $n$-dimensional isotropic subspaces of the Cartan distribution of a PDE. These are tangent to solutions of PDEs. In this way we are able to generalize ``Maslov index" for Lagrangian submanifolds as introduced by V. I. Arnold, to any solution of PDEs. Really Lagrangian submanifolds of symplectic manifolds, can be encoded as solutions of suitable first order PDEs.

As a by-product we get also a new proof for existence of the Navier-Stokes PDEs global smooth solutions, defined on all $\mathbb{R}^3$. (Example \ref{navier-stokes-pdes-and-global-space-time-smooth-solutions}.)

In this section we shall calculate also Lagrangian bordism groups in a $2n$-dimensional symplectic manifold $(W,\omega)$, where $\omega$ is a non-degenerate, close, differentiable $2$-differential form on $W$. In \cite{PRA022} we have calculated the Lagrangian bordism groups in the case that $\omega$ is exact. This has been made by generalizing to higher order PDE, a previous approach given by V. I. Arnold \cite{ARNOLD1,ARNOLD2}, Y. Eliashberg \cite{ELIASHBERG} and A. Pr\'astaro \cite{PRA05}. Now we give completely new formulas, without assuming any restriction on $\omega$, and following our Algebraic Topology of PDEs. (See References \cite{PRA05,PRA022,PRA00,PRA000,PRA01,PRA1,PRA022,PRA3,PRA4,PRA5,PRA6,PRA7,PRA8,PRA9,PRA10,PRA11,PRA14,PRA15,PRA16,PRA17,PRA017}. See also \cite{AG-PRA1,AG-PRA2,LYCH-PRAS,PRA-RAS}.)

In this section our main results are Theorem \ref{maslov-index-and-maslov-cycle-relation-solution-pde}, Theorem \ref{maslov-index-for-lagrangian-manifolds}, Theorem \ref{g-singular-lagrangian-bordism-groups} and Theorem \ref{weak-and-singular-lagrangian-bordism-groups-a}. The first is devoted to relation between Maslov indexes and Maslov cycles for solutions of PDEs. The second characteriizes such invariants for Lagrangian submanifolds of symplectic manifolds, by means of suitable formally integrable and completely integrable first order PDEs. The other two theorems characterize Lagrangian bordism groups in such PDEs.

\begin{theorem}[Grassmannian of $n$-dimensional integral planes of $J^k_n(W)$]\label{grassmannian-of-n-dimensional-integral-planes}
Let $I_k(W)_q$ be the {\em Grassmannian of $n$-dimensional integral planes}  at $q\in J^k_n(W)$, namely the set of isotropic $n$-dimensional subspaces of the Cartan distribution $\mathbf{E}_k(W)_q$. One has the following properties.

{\em(i)} One has the natural fiber bundle structure $ I_k(W)=\bigcup_{q\in J^k_n(W)}I_k(W)_q\to J^k_n(W)$.

{\em(ii)} In general an integral $n$-plane $L\in I_k(W)_q$, is projected, via $(\pi_{k,0})_*$ onto an $(n-l)$-dimensional subspace of $T_aN$, $q=[N]^k_a$.

{\em(iii)} The set of of $n$-integral planes such that $\dim(\pi_{k,0})_*(L)=n=\dim(T_aN)$, (namely with $l=0$), is identified with the affine fiber $\pi^{-1}_{k+1,k}(q)\subset J^{k+1}_n(W)$. These integral planes are called {\em regular integral planes}.

{\em(iv)}  In general an $n$-integral plane $L\in  I_k(W)_q$, admits the following splitting
\begin{equation}\label{splitting-n-integral-plane}
  L\cong L_o\bigoplus L_v
\end{equation}
where $L_o$ {\em(horizontal component)}, is contained in some regular plane $L_{\widetilde{q}}$, for some $\widetilde{q}\in\pi^{-1}_{k+1,k}(q)\subset J^{k+1}_n(W)$. Furthermore $L_v$, {\em(vertical component)}, is contained in the vector space $T_q\pi^{-1}_{k,k-1}(\bar q)\cong S^k(T_a^*N)\bigotimes \nu_a$, with $\bar q=\pi_{k,k-1}(q)\in J^{k-1}_n(W)$.

{\em(v)}  Two different splittings $L\cong L_o\bigoplus L_v$ and $L\cong L'_o\bigoplus L'_v$, of a $n$-integral plane $L\subset \mathbf{E}_k(W)_q$, $q=[N]_a^k\in J^k_n(W)$, are related by a fixed subspace $V\subset S^k(T^*_aN)\bigotimes\nu_a$. More precisely one has:
\begin{equation}\label{relation-between-two-splittings}
  L_o=L'_o\bigoplus V;\, L'_v=L_v\bigoplus V.
\end{equation}

$\bullet$\hskip 2pt {\em(Cohomology ring $H^\bullet(I_k(W))$)}. One has the following isomorphisms:
\begin{equation}\label{cohomology-ring-formulas-isomorphisms-integral-planes}
  \left\{
  \begin{array}{ll}
  H^\bullet(I_k(W);\mathbb{Z}_2)&\cong H^\bullet(J^k_n(W);\mathbb{Z}_2)\otimes_{\mathbb{Z}_2}H^\bullet(F_k(W);\mathbb{Z}_2)\\
  &\cong H^\bullet(W;\mathbb{Z}_2)\otimes_{\mathbb{Z}_2}H^\bullet(F_k(W);\mathbb{Z}_2)\\
  \end{array}
  \right.
\end{equation}
where $F_k(W)$ is the fiber of $I_k(W)$ over $J^k_n(W)$. One has the following ring isomorphism:
$$H^\bullet(F_k(W);\mathbb{Z}_2)\cong\mathbb{Z}_2[w^{(k)}_1,\cdots,w^{(k)}_n],$$
where $\deg(w^{(k)}_i)=i$. Such generators coincide with Stiefel-Whitney classes of the tautological bundle $E(\eta)\to I_k(W)$.
\end{theorem}

\begin{proof}
Let us only explicitly consider that the first part of the formula (\ref{cohomology-ring-formulas-isomorphisms-integral-planes}) follows from a direct application of some results about spectral sequences an their relations with fibration ({\em Leray-Hirsh theorem}). For more details see Theorem 3 in \cite{PRA00}.
\end{proof}

\begin{theorem}[Grassmannian of $n$-dimensional integral planes of PDE]\label{grassmannian-of-n-dimensional-integral-planes-of-pde}
Let
$$I(E_k)=\bigcup_{q\in E_k}I(E_k)_q$$
be the Grassmannian of $n$-dimensional integral planes of $E_k$. One has a natural fiber bundle structure $I(E_k)\to E_k$. Then each singular solution $V\subset E_k$ identifies a mapping $i_V:V\to I(E_k)$, given by $i_V(q)=T_qV\in I(E_k)_q$. Then one has an induced morphism
\begin{equation}\label{induced-morphism}
  i_V^*:H^i(I(E_k):\mathbb{Z}_2)\to H^i(V;\mathbb{Z}_2),\, \omega\mapsto i_V^*\omega.
\end{equation}
$i_V^*\omega$ is the {\em characteristic class} of $V$ corresponding to $\omega$.

$\bullet$\hskip 2pt If $E_k$ is a strong retract of $J^k_n(W)$ then $H^\bullet(I(E_k);\mathbb{Z}_2)$ is an algebra over $H^\bullet(E_k;\mathbb{Z}_2)$.
More precisely one has $$ H^i(I(E_k);\mathbb{Z}_2)\cong \bigoplus_{r+s=i}H^r(E_k;\mathbb{Z}_2)\bigotimes_{\mathbb{Z}_2}H^s(F_k;\mathbb{Z}_2),$$ where
$F_k$ is the fibre of $I(E_k)$ over $E_k$.

$\bullet$\hskip 2pt Furthermore, the ring $H^\bullet(F_k;\mathbb{Z}_2)$ is isomorphic up to $n$ to the ring $\mathbb{Z}_2[\omega^{(k)}_1,\cdots,\omega^{(k)}_n)$ of polynomials in the generator $\omega^{(k)}_i$, ${\rm degree}(\omega^{(k)}_i)=i$. These generators can be identified with the Stiefel-Whitney classes of the tautological bundle $E(\eta)\to I(E_k)$.

$\bullet$\hskip 2pt If $E_k$ is a formally integrable PDE then

\begin{equation}\label{formally-integrability-cohomology-integral-planes-pde}
\left\{\begin{array}{ll}
H^\bullet(I(E_{k+1});\mathbb{Z}_2)&\cong H^\bullet(I_{k+1}(W);\mathbb{Z}_2)\\
&\cong H^\bullet(W;\mathbb{Z}_2))\bigotimes_{\mathbb{Z}_2}H^\bullet(F_{k+1}(W);\mathbb{Z}_2)\\
&\cong H^\bullet(W;\mathbb{Z}_2))\bigotimes_{\mathbb{Z}_2}\mathbb{Z}_2[\omega^{(k+1)}_1,\cdots,\omega^{(k+1)}_n).\\
\end{array}\right.
\end{equation}

$\bullet$\hskip 2pt If $V$ is a non-singular solution of $E_k$, then all its characteristic classes are zero in dimension $\ge 1$.
\end{theorem}
\begin{proof}
After Theorem \ref{grassmannian-of-n-dimensional-integral-planes}, let us only explicitly consider when $E_k$ is a strong retract of $J^k_n(W)$.
This fact implies the homotopy equivalence $E_k\backsimeq J^k_n(W)$. Then we can state also the homotopy equivalence between the corresponding integral planes fiber-bundles $I(E_k)\backsimeq I_k(W)$. In fact we use the following lemmas.

\begin{lemma}\label{homotopy-equivalence-strong-retract}
If $A\subset X$ is a strong retract, then the inclusion $i:(A,x_0)\hookrightarrow (X,x_0)$ is an homotopy equivalence and hence $i_*:\pi_n(A,x_0)\to\pi_n(X,x_0)$ is an isomorphism for all $n\ge 0$.
\end{lemma}
\begin{proof}
This is a standard result. See e.g., \cite{PRA3}. (This lemma is the inverse of the {\em Whitehead's theorem}.)
\end{proof}
\begin{lemma}\label{homotopy-equivalence-fiber-bundles-from-homotopy-equivalence-between-bases}
For a space $B$ let $\mathcal{F}(B)$ be the set of fiber homotopy equivalence classes of fibrations $E\to B$. A map $f:B_1\to B_2$ induces $f^*:\mathcal{F}(B_2)\to\mathcal{F}(B_1)$, depending only on the homotopy class of $f$. If $f$ is a homotopy equivalence, then $f^*$ becomes a bijection: $f^*:\mathcal{F}(B_2)\leftrightarrow\mathcal{F}(B_1)$.
\end{lemma}
\begin{proof}
This is a standard result. See, e.g., \cite{HATCHER}.
\end{proof}

From above two lemmas, we can state that also $I(E_k)$ is a strong retract of $I_k(W)$, therefore one has the following exact commutative diagram of homotopy equivalences:

\begin{equation}\label{homotopy-equivalence-comm-diagram-strong-retract-pde-integral-planes}
  \xymatrix@C=2cm{0\ar[r]&I(E_k)\ar[d]\ar@{-}[r]^{\sim}&I_k(W)\ar[d]\\
  0\ar[r]&E_k\ar[d]\ar@{-}[r]^{\sim}&J_n^k(W)\ar[d]\\
  &0&0\\}
\end{equation}
This induces the following commutative diagram of isomorphic cohomologies:

\begin{equation}\label{homotopy-equivalence-comm-diagram-strong-retract-pde-integral-planes-a}
  \xymatrix@C=2cm{H^\bullet(I(E_k);\mathbb{Z}_2)\ar@{=}[d]^{\wr}\ar@{=}[r]^{\sim}&H^\bullet(I_k(W);\mathbb{Z}_2)\ar@{=}[d]^{\wr}\\
  H^\bullet(E_k;\mathbb{Z}_2)\bigotimes_{\mathbb{Z}_2}H^\bullet(F_k;\mathbb{Z}_2)
  \ar@{=}[r]^{\sim}&H^\bullet(J_n^k(W);\mathbb{Z}_2)\bigotimes_{\mathbb{Z}_2}H^\bullet(F_k(W);\mathbb{Z}_2)\\}
\end{equation}
Since $$H^\bullet(E_k;\mathbb{Z}_2)\cong H^\bullet(J^k_n(W);\mathbb{Z}_2)\cong H^\bullet(W;\mathbb{Z}_2)$$ and
$$H^\bullet(F_k;\mathbb{Z}_2)\cong H^\bullet(F_k(W);\mathbb{Z}_2)\cong \mathbb{Z}_2[\omega_1^{(k)},\cdots,\omega_n^{(k)}],$$ we get
$$H^\bullet(I(E_k);\mathbb{Z}_2)\cong H^\bullet(W;\mathbb{Z}_2)\bigotimes_{\mathbb{Z}_2}\mathbb{Z}_2[\omega_1^{(k)},\cdots,\omega_n^{(k)}].$$
Therefore, $H^\bullet(I(E_k);\mathbb{Z}_2)$ is an algebra over $H^\bullet(W;\mathbb{Z}_2)$ isomorphic to $\mathbb{Z}_2[\omega_1^{(k)},\cdots,\omega_n^{(k)}]$.

Finally if $E_k$ is formally integrable, then its $r$-prolongations $(E_k)_{+r}$ are strong retract of $J^{k+r}_n(W)$, for $r\ge 1$. Thus we can repeat above considerations by working on each $(E_k)_{+r}$ and obtain
$$H^\bullet(I((E_k)_{+r});\mathbb{Z}_2)\cong H^\bullet(W;\mathbb{Z}_2)\bigotimes_{\mathbb{Z}_2}\mathbb{Z}_2[\omega_1^{(k+r)},\cdots,\omega_n^{(k+r)}],$$ for $r\ge 1$.
\end{proof}

\begin{remark}\label{comparison-between-metasymplectic-structure-symplectic-structure}
It is worth to emphasize the comparison between metasymplectic strucure on $J^k_n(W)$, and the symplectic structure in a symplectic vector space $(V,\omega)$. According to the definition given in the proof of Theorem \ref{metasymplectic-structure}, we can define {\em metasymplectic orthogonal} of a subspace $P\triangleleft\mathbf{E}_k(W)_q$, the set

\begin{equation}\label{methasymplectic-orthogonal}
\left\{
\begin{array}{ll}
P^{\bot}&=\{\zeta\in\mathbf{E}_k(W)_q\, |\, \Omega_k(\lambda)(\zeta,\xi)=0,\, \forall \xi\in P,\, \forall\lambda\in S^{k-1}(T_aN)\otimes\nu_a^*\}\\
&=\bigcap_{\lambda\in S^{k-1}(T_aN)\otimes\nu_a^*}\ker(\Omega_k(\lambda)(\zeta,P)).\\
\end{array}
\right.
\end{equation}

One has the following properties:

{\em(a)} $(P^{\bot})^{\bot}=P$;

{\em(b)} $P_1^{\bot}\bigcap P_2^{\bot}=(P_1+P_2)^{\bot}$;

{\em(c)} $(P_1\bigcap P_2)^{\bot}=(P_1)^{\bot}+(P_2)^{\bot}$.

Then one can define $P$ {\em metasymplectic-isotropic} if $P\subset P^{\bot}$. Furthermore, we say that $P$ is {\em metasymplectic-Lagrangian} if $P=P^{\bot}$. Maximal metasymplectic-isotropic spaces are metasymplectic-isotropic spaces that are not contained into larger ones. There any two vectors are an involutive couple. With respect to above remarks, in Tab. \ref{table-comparison-between-metasymplectic-structure-symplectic-structure} we have made a comparison between definitions related to the metasymplectic structure and symplectic structure. Let us underline that the metasymplectic structure considered, is not a trivial extension of the canonical symplectic structure that can be recognized on any vector space $E$, of dimension $n$. In fact, it is well known that $V=E\bigoplus E^*$, has the canonical symplectic structure $\sigma((v,\alpha),(v',\alpha'))=<\alpha,v'>-<\alpha',v>$, called the {\em natural symplectic form} on $E$. Instead the metasymplectic structure arises from differential of Cartan forms.
\end{remark}

\begin{table}[t]
\caption{Comparison between metasymplectic structure of $J^k_n(W)$ and symplectic structure of symplectic space $(V,\omega)$.}
\label{table-comparison-between-metasymplectic-structure-symplectic-structure}
\scalebox{0.8}{$\begin{tabular}{|l|l|}
\hline
\hfil{\rm{\footnotesize $(\mathbf{E}_k(W)_q\, ,\, \Omega_k(\lambda))$}}\hfil&\hfil{\rm{\footnotesize $(V,\omega)$}}\\
\hfil{\rm{\footnotesize $\dim (\mathbf{E}_k(W)_q)=m\binom{k+n-1}{k}+n$, $\dim W=n+m$}}\hfil&\hfil{\rm{\footnotesize $\dim V=2n$}}\hfil\\
\hline\hline
\hfil{\rm{\footnotesize $P\triangleleft\mathbf{E}_k(W)_q$}}\hfil&\hfil{\rm{\footnotesize $E\triangleleft V$}}\hfil\\
\hfil{\rm{\footnotesize $P^{\bot}=\{\zeta\in \triangleleft\mathbf{E}_k(W)_q\, |\, \Omega_k(\lambda)(\zeta,\xi)=0,\, \forall \xi\in P,\, \forall\lambda\in S^{k-1}(T_aN)\otimes\nu_a^*\}$}}\hfil&\hfil{\rm{\footnotesize $E^{\bot}=\{v\in V\, |\, \omega(v,u)=0,\, \forall u\in E\}$}}\hfil\\
\hline
{\rm{\footnotesize $P$ metasymplectic-isotropic iff $P\subset P^{\bot}$.}}&{\rm{\footnotesize $E$ is symplectic-isotropic iff $E\subset E^{\bot}$.}}\\
\hfil{\rm{\footnotesize $[\dim P\le P^{\bot}]$}}\hfil&\hfil{\rm{\footnotesize $[\dim E\le n]$}}\hfil\\
\hline
{\rm{\footnotesize $P$ is metasymplectic-Lagrangian iff $P=P^{\bot}$.}}&{\rm{\footnotesize $E$ is symplectic-Lagrangian iff $E=E^{\bot}$.}}\\
\hfil{\rm{\footnotesize $[\dim P=P^{\bot}]$}}\hfil&\hfil{\rm{\footnotesize $[\dim E= n]$}}\\
\hline
{\rm{\footnotesize $P$ is maximal metasymplectic-isotropic iff  $P \not\subset Q\subset Q^{\bot}$.}}&{\rm{\footnotesize $E$ is symplectic-co-isotropic iff $E\supset E^{\bot}$.}}\\
\hfil{\rm{\footnotesize  $[\dim P=m\binom{p+k-1}{k}+n-p$, $0\le p\le n]$}}\hfil&\hfil{\rm{\footnotesize $[\dim E\ge n]$}}\hfil\\
\hfil{\rm{\footnotesize  $[\hbox{\rm type} (P)=n-p]$}}\hfil&\hfil{\rm{\footnotesize }}\hfil\\
\hline
\multicolumn{2}{l}{\rm{\footnotesize A metasymplectic-isotropic space is metasymplectic-involutive.}}\hfill\\
\multicolumn{2}{l}{\rm{\footnotesize A metasymplectic-Lagrangian space is metasymplectic-involutive.}}\hfill\\
\multicolumn{2}{l}{\rm{\footnotesize A symplectic-Lagrangian space is maximally symplectic-isotropic.}}\hfill\\
\multicolumn{2}{l}{\rm{\footnotesize A symplectic-isotropic (or symplectic-co-isotropic) space $E$ with $\dim E=n$, is symplectic-Lagrangian.}}\hfill\\
\multicolumn{2}{l}{\rm{\footnotesize A line (hyperplane) is symplectic-isotropic (symplectic-co-isotropic).}}\hfill\\
\multicolumn{2}{l}{\rm{\footnotesize A maximal metasymplectic-isotropic space of type $n$ has dimension $n$.}}\hfill\\
\multicolumn{2}{l}{\rm{\footnotesize A maximal metasymplectic-isotropic space of type $0$ has dimension $m\binom{n+k-1}{k}=\dim[S^k(T^*_aN)\otimes\nu_a]$.}}\hfill\\
\end{tabular}$}
\end{table}

\begin{definition}[Lagrangian submanifolds of symplectic manifold]\label{lagrangian-submanifolds-of-symplectic-manifold}
Let $(W,\omega)$ be a {\em symplectic manifold}, that is $W$ is a $2n$-dimensional manifold with symplectic $2$-form $\omega:W\to \Lambda^0_2(W)$, (hence $\omega$ is closed: $d\omega=0$). We call Lagrangian manifold a $n$-dimensional submanifold $V\subset W$, such that $\omega|_V=0$.\footnote{The tangent space $T_pW$, $\forall p\in W$, identifies a symplectic space via the $2$-form $\omega(p)\in \Lambda^2(T^*_pW)$. Therefore a $n$-dimensional sub-manifold $V$ of a $2n$-symplectic manifold $W$ is Lagrangian iff $T_pV$ is a Lagrangian subspace of $T_pW$, $\forall p\in V$.}
\end{definition}

\begin{example}{\em(Arnold 1967)}\label{example-lagrangian-submanifold-a}
A Lagrangian submanifold of the symplectic space $(\mathbb{R}^{2n},\omega)$ is a $n$-dimensional submanifold $V\subset \mathbb{R}^{2n}$, such that for any $p\in V$, $T_pV\subset T_p\mathbb{R}^{2n}\cong \mathbb{R}^{2n}$, is a Lagrangian subspace of $\mathbb{R}^{2n}$. This is equivalent to say that the symplectic $2$-form $\sigma=\sum_{1\le r<s\le 2n}\sigma_{rs}d\xi^r\wedge d\xi^s$, with $\sigma_{rs}=\omega_{rs}$, and $(\xi^r)_{1\le r\le 2n}=(x^j,y^j)_{1\le j\le n}$, annihilates on $V$: $\sigma|_V=0$. The tangent space $TV$, classified by the first classifying mapping $f:V\to BO(n)$, is the pullback of the tautological bundle $E(\eta)$ over $L_{agr}(\mathbb{R}^{2n},\omega) $, or equivalently the pullback of $E(\eta)$ via the second classifying mapping $\zeta:V\to L_{agr}(\mathbb{R}^{2n},\omega)$, $\zeta(p)=T_pV\cong\mathbb{R}^{2n}$. In fact one has the exact commutative diagram {\em(\ref{exact-commutative-diagram-example-lagrange-submanifold})}.
\begin{equation}\label{exact-commutative-diagram-example-lagrange-submanifold}
  \xymatrix{TM\cong f^* E(\eta)\cong \zeta^* E(\eta)\ar[d]\ar[r]\ar@/ ^2pc/[rr]&E(\eta)\ar[d]\ar@{=}[r]&E(\eta)\ar[d]\\
  V\ar[d]\ar[r]_{\zeta}\ar@/ _2pc/[rr]_{f}&L_{agr}(\mathbb{R}^{2n},\omega)\ar[d]\ar[r]_{\eta}&BO(n)\ar[d]\\
  0&0&0\\}
\end{equation}

$\bullet$\hskip 2pt The {\em Maslov index class} of $V$ is defined by $\tau(V)=\zeta^*(\alpha)\in H^1(V;\mathbb{Z})$, where $\alpha\in H^1(L_{agr}(\mathbb{R}^{2n},\omega);\mathbb{Z})\cong\mathbb{Z}$ is the generator.

$\bullet$\hskip 2pt The {\em Maslov cycle} of $V$ is defined by
$$\Sigma(V)=\{p\in V \, |\, \dim(\ker(T(\pi)|_{T_pV}))>0\}$$
where $\pi:\mathbb{R}^{2n}\cong\mathbb{R}^{n}\bigoplus i\mathbb{R}^{n}\to\mathbb{R}^{n}$. Therefore $\Sigma(V)\cong T_pV\bigcap i\mathbb{R}^{n}\not=\{0\}$. The homology class $[\Sigma(V)]\in H_{n-1}(V;\mathbb{Z})$ is the Poincar\'e dual of the  Maslov index class $\tau(V)\in H^{1}(V;\mathbb{Z})$:
$$[\Sigma(V)]=D\tau(V).$$

Therefore, one can state that $\tau(V)$ measures the failure of the morphism $\pi|_V:V\to \mathbb{R}^n$ to be a local diffeomorphism.
\end{example}
\begin{example}\label{example-lagrangian-submanifold-b}
$\mathbb{C}^n$ is a symplectic manifold. Any $n$-dimensional subspace is a Lagrangian submanifold.
\end{example}

\begin{example}\label{example-lagrangian-submanifold-c}
Any $1$-dimensional submanifold of a $2$-dimensional symplectic manifold is Lagrangian.\footnote{For example any curve in $S^2$ is a Lagrangian submanifold.}
\end{example}

\begin{example}\label{example-lagrangian-submanifold-d}
The cotangent space $T^*M$ of a $n$-dimensional manifold $M$ is a symplectic manifold, and each fiber $T_p^*M$ of the fiber bubdle $\pi:T^*M\to M$, is a Lagrangian submanifold.

$\bullet$\hskip 2pt Let $V\subset T^*M$ be a Lagrangian submanifold of $T^*M$. Let us consider the fiber bundle
\begin{equation}\label{L_{agr}angian-bundle-cotangent-bundle}
  L_{agr}(T^*M)=\bigcup_{q\in T^*M}L_{agr}(T^*M)_q
\end{equation}
where $L_{agr}(T^*M)_q$ is the set of Lagrangian subspaces of $T_q(T^*M)$. One has a canonical mapping
$$\zeta:V\to L_{agr}(T^*M),\, q\mapsto T_qV.$$
Then if $\alpha\in H^1(T^*M;\mathbb{Z})\cong\mathbb{Z}$ is the generator, we get $\zeta^*\alpha\in H^1(V;\mathbb{Z})$ is the {\em Maslov index class} of $V$. The {\em Maslov cycle} of $V$ is defined the set
$$\Sigma(V)=\{q\in V\, |\, \dim (\ker(T(\pi)|_{T_qV}>0, \, \pi:T^*M\to M\}.$$
Therefore $\Sigma(V)\cong\{q\in V\, |\, T_qV\bigcap vT_q(T^*M)\not=\{0\}$. Here $vT_q(T^*M)$ denotes the vertical tangent space at $q\in T^*M$, with respect to the projection $\pi:T^*M\to M$. The homology class $[\Sigma(V)]\in H_{n-1}(V;\mathbb{Z})$ is the Poincar\'e dual of the Maslov index class $\tau(V)\in H^1(V;\mathbb{Z})$. $[\Sigma(V)]=D\tau(V)$. Therefore $\tau(V)$ measures the failure of the mapping $\pi|_V:V\to M$ to be a local diffeomorphism.
\end{example}

\begin{definition}[Maslov cycles of PDE solution]\label{maslov-cycle-solution-pde}
We call {\em $i$-Maslov cycle}, $1\le i\le n-1$, of a solution $V\subset E_k\subset J^k_n(W)$, the set $\Sigma_i(V)$ of singular points $q\in V$, such that $\dim(\ker((\pi_{k,0})_*)|_{T_qV})))=n-i$.
\end{definition}

\begin{definition}[Maslov index classes of PDE solution]\label{maslov-index-class-solution-pde}
We call {\em $i$-Maslov index class}, $1\le i\le n-1$, of a solution $V\subset E_k\subset J^k_n(W)$,  $$\tau_i(V)=(i_V)^*\omega_i^{(k)}\in H^i(V;\mathbb{Z}_2),$$ where $\omega_i^{(k)}$ is the $i$-th generators of the ring $\mathbb{Z}_2[\omega_1^{(k)},\cdots,\omega_n^{(k)}]\cong H^\bullet(F_k,\mathbb{Z}_2)$ and $i_V:V\to I(E_k)$ is the canonical mapping, $i_V:q\mapsto T_qV$.
\end{definition}
\begin{theorem}[Maslov indexes and Maslov cycles relations for solution of PDE]\label{maslov-index-and-maslov-cycle-relation-solution-pde}

$\bullet$\hskip 2pt Let $E_k\subset J^k_n(W)$ be a strong retract of $J^k_n(W)$, then
the homology class $[\Sigma_i(V)]\in H_{n-i}(V;\mathbb{Z})$, $1\le i\le n-1$, is the Poincar\'e dual of the Maslov index class $\tau_i(V)\in H^i(V;\mathbb{Z})$. Formula {\em(\ref{formula-dpoincare-duality-main})} holds.
\begin{equation}\label{formula-dpoincare-duality-main}
  [\Sigma_i(V)]=D\tau_i(V),\, 1\le i\le n-1.
\end{equation}
Therefore, $\{\tau_i(V)\}_{1\le i\le n-1}$, measure the failure of the mapping $\pi_{k,0}:V\to W$ to be a local embedding.

$\bullet$\hskip 2pt Let $E_k\subset J^k_n(W)$ be a formally integrable PDE. Then one can charcaterize each solution $V$ on the first prolongations $(E_k)_{+1}\subset J^{k+1}_n(W)$, by means of $i$-Maslov indexes and $i$-Maslov cycles, as made in above point.
\end{theorem}

\begin{proof}
Let us consider $E_k$ a strong retract of $J^k_n(W)$. Then we can apply Theorem \ref{grassmannian-of-n-dimensional-integral-planes-of-pde}
In particular we get the following isomorphisms:
\begin{equation}\label{cohomology-isomorphisms-integral-planes}
\left\{
 \begin{array}{ll}
 H^\bullet(E_k)&\cong H^\bullet(J^k_n(W))\\
 &\cong H^\bullet(W)\\
 H^\bullet(I(E_k))&\cong H^\bullet(I_k(W)).\\
  \end{array}
 \right.
 \end{equation}
Let us more explicitly calculate these cohomologies. Start with the case $i=1$. One has the following isomorphisms:
\begin{equation}\label{cohomology-isomorphisms-integral-planes-a}
\left\{
 \begin{array}{ll}
  H^1(I(E_k);\mathbb{Z}_2)&\cong H^1(E_k;\mathbb{Z}_2)\otimes_{\mathbb{Z}_2}H^0(F_k;\mathbb{Z}_2)\bigoplus H^0(E_k;\mathbb{Z})\otimes_{\mathbb{Z}_2}H^1(F_k;\mathbb{Z}_2)\\
  &\cong H^1(E_k;\mathbb{Z}_2)\otimes_{\mathbb{Z}_2}\mathbb{Z}_2\bigoplus \mathbb{Z}_2\otimes_{\mathbb{Z}_2}\mathbb{Z}_2[\omega_1^{(k)}]\\
  &\cong H^1(E_k;\mathbb{Z}_2)\bigoplus \mathbb{Z}_2[\omega_1^{(k)}].\\
  \end{array}
 \right.
 \end{equation}
Therefore one has the following exact commutative diagram:
 \begin{equation}\label{commutative-diagram-cohomologies-integral-planes-a}
 \xymatrix{0\ar[r]&\mathbb{Z}_2[\omega_1^{(k)}]\ar[dr]_{(i_V)_*=(i_V)^*|_{\mathbb{Z}_2[\omega_1^{(k)}]}}\ar[r]&H^1(I(E_k);\mathbb{Z}_2)\ar[d]_{(i_V)^*}\ar[r]&H^1(E_k;\mathbb{Z}_2)\ar[r]&0\\
 &&H^1(V;\mathbb{Z}_2)\ar[d]&&\\
 &&0&&\\}
 \end{equation}
Then the mapping $i_V:V\to I(E_k)$, induces the following morphism $$(i_V)_*=(i_V)^*|_{\mathbb{Z}_2[\omega_1^{(k)}]}:\mathbb{Z}_2[\omega_1^{(k)}]\to H^1(V;\mathbb{Z}_2).$$
Set $\beta_1(V)=(i_V)_*(\omega_1^{(k)})$. Here we suppose that $V$ is compact, (otherwise we shall consider cohomology with compact support). Now we get
\begin{equation}\label{proof-a}
\beta_1(V)\bigcap[V]=[\Sigma_1(V)].
\end{equation}
In (\ref{proof-a}) $[V]$ denotes the fundamental class of $V$ that there exists also whether $V$ is non-orientable. (For details see, e.g. \cite{PRA3}.)

We can pass to any degree, $1\le i\le n-1$, by considering the following isomorphisms:

\begin{equation}\label{cohomology-isomorphisms-integral-planes-b}
\left\{
 \begin{array}{ll}
  H^i(I(E_k);\mathbb{Z}_2)&\cong H^i(E_k;\mathbb{Z}_2)\\
  &\bigoplus_{1\le p\le i-1}H^{i-p}(E_k;\mathbb{Z}_2)\otimes_{\mathbb{Z}}H^p(F_k;\mathbb{Z}_2)\\
  &\bigoplus \mathbb{Z}_2[\omega_1^{(k)},\cdots,\omega_{i}^{(k)}].\\
  \end{array}
 \right.
\end{equation}
 One has the following exact commutative diagram:

 \begin{equation}\label{commutative-diagram-cohomologies-integral-planes-b}
  \scalebox{0.8}{$\xymatrix{0\ar[r]&\mathbb{Z}_2[\omega_1^{(k)},\cdots,\omega_i^{(k)}]\ar[dr]_{(i_V)_*}\ar[r]&H^i(I(E_k);\mathbb{Z}_2)\ar[d]_{(i_V)^*}\ar[r]&
  \framebox{$H^i(E_k;\mathbb{Z}_2)/\mathbb{Z}_2[\omega_1^{(k)},\cdots,\omega_i^{(k)}]$}\ar[r]&0\\
 &&H^i(V;\mathbb{Z}_2)\ar[d]&&\\
 &&0&&\\}$}
 \end{equation}
 Then the map $i_V:V\to I(E_k)$ induces  the following morphism:
 $$(i_V)_*:\mathbb{Z}_2[\omega_1^{(k)},\cdots,\omega_i^{(k)}]\to H^i(V;\mathbb{Z}_2),\, 1\le i\le n-1.$$
 Set $\beta_i(V)=(i_V)_*(\omega_i^{(k)})$. We get
\begin{equation}\label{proof-a}
\beta_i(V)\bigcap[V]=[\Sigma_i(V)].
\end{equation}

For the case where $E_k$ is formally integrable, we can repeat the above proceeding applied to the first prolongation $(E_k)_{+1}$ of $E_k$, that is a strong retract of $J^{k+1}_n(W)$. In this way we complete the proof.
\end{proof}

\begin{figure}[t]
\centerline{\includegraphics[width=8cm]{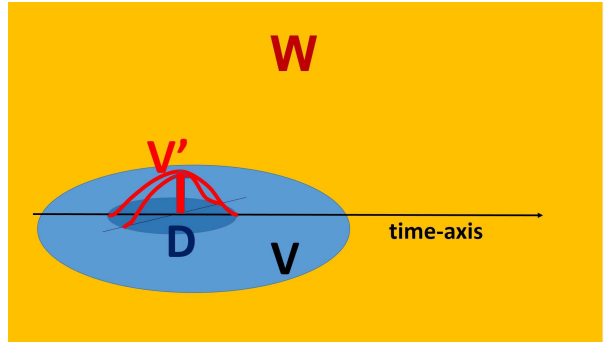}}
\renewcommand{\figurename}{Fig.}
\caption{Global space-time smooth solution representation $V'\subset \widehat{(NS)}$, obtained by means of a localized, space-time smooth perturbation of a constant global smooth solution $V\subset\widehat{(NS)}$. The perturbation, localized in the compact space-time region $D\subset V$ of the smooth global constant solution $V$, is a smooth solution of the linearized equation $\widehat{\overline{(NS)}}[s]$. The vertical arrow denotes the local perturbation of the solution $V$, generating the non-constant global smooth solution $V'\subset J^2_4(W)$.}
\label{global-nonconstant-smooth-solution-representation}
\end{figure}

\begin{example}[Navier-Stokes PDEs and global space-time smooth solutions]\label{navier-stokes-pdes-and-global-space-time-smooth-solutions}
The non-isothermal Navier-Stokes equation can be encoded in a geometric way as a second-order PDE $(NS)\subset J^2_4(W)$, where $\pi:W= J{\it D}(M)\times_M T^0_0M\times_M T^0_0M\cong M\times \mathbf{I}\times\mathbb{R}^2\to M$ is an affine bundle over the $4$-dimensional affine Galilean space-time $M$. There $\mathbf{I}\subset \mathbf{M}$ represents a $3$-dimensional affine subspace of the $4$-dimensional vector space $\mathbf{M}$ of free vectors of $M$. A section $s:M\to W$ is a triplet $s=(v,p,\theta)$ representing the velocity field $v$, the isotropic pressure $p$, and the temperature $\theta$. In \cite{PRA1} it is reported the explicit expression of $(NS)$, formulated just in this geometric way. Then one can see there that $(NS)$ is not formally integrable, but one can canonically recognize a sub-equation $\widehat{(NS)}\subset (NS)\subset J^2_4(W)$, that is so and also completely integrable. Furthermore, $(NS)$ is a strong deformed retract of $J^2_4(W)$, over a strong deformed retract $(C)$ of $J^1_4(W)$.  In other words one has the following commutative diagram of homotopy equivalences:
\begin{equation}\label{homotopy-equivalences-navier-stokes-pdes}
  \xymatrix@R=1cm @C=2cm{(NS)\ar[d]\ar@{-}[r]^{\sim}&J^2_4(W)\ar[d]\\
  (C)\ar[d]\ar@{-}[r]^{\sim}&J^1_4(W)\ar[d]\\
  0&0\\}
\end{equation}
Since $J_4^2(W)$ and $J^1_4(W)$ are affine spaces, they are topologically contractible to a point, hence from {\em(\ref{homotopy-equivalences-navier-stokes-pdes})} we are able to calculate the cohomology properties of $(NS)$, as reported in {\em(\ref{cohomology-space-navier-stokes-pdes})}.
\begin{equation}\label{cohomology-space-navier-stokes-pdes}
 \left\{
 \begin{array}{l}
   H^0((NS);\mathbb{Z}_2)=\mathbb{Z}_2\\
  H^r((NS);\mathbb{Z}_2)=0,\, r>0.\\
  \end{array}
 \right.
\end{equation}
We get the cohomologies of $I(NS)$, as reported in {\em(\ref{cohomology-space-navier-stokes-pdes-a})}.
\begin{equation}\label{cohomology-space-navier-stokes-pdes-a}
 \left\{
 \begin{array}{ll}
  H^r(I(NS);\mathbb{Z}_2)&=\bigoplus_{p+q=r}H^p((NS);\mathbb{Z}_2)\otimes_{\mathbb{Z}_2}H^q(F_2;\mathbb{Z}_2)\\
  &=H^0((NS);\mathbb{Z}_2)\otimes_{\mathbb{Z}_2}H^r(F_2;\mathbb{Z}_2)\\
  &=\mathbb{Z}_2\otimes_{\mathbb{Z}_2} \mathbb{Z}_2[\omega_{1}^{(2)},\cdots \omega_{r}^{(2)}]\\
  &=\mathbb{Z}_2[\omega_{1}^{(2)},\cdots \omega_{r}^{(2)}],\, 1\le r\le 4.
  \end{array}
 \right.
\end{equation}

Therefore {\em(\ref{cohomology-space-navier-stokes-pdes-c})} are the conditions that $V\subset \widehat{(NS)}\subset(NS)$ must satisfy in order to be without singular points.
\begin{equation}\label{cohomology-space-navier-stokes-pdes-c}
  0=i_V^*\omega_i^{(2)}\in H^i(V;\mathbb{Z}_2),\, 1\le i\le 4.
\end{equation}
In particular, if $$V=D^2s(M)\subset \widehat{(NS)}\subset (NS)\subset J{\it D}^2(W)\subset J^2_4(W)$$
where $s:M\to W$ is a smooth global section, since $H^i(V;\mathbb{Z}_2)=0$, $\forall i>0$,  we get that all its characeristic classes $i_V^*\omega_i^{(2)}$ are zero. Therefore, $V$ cannot have singular points on $V$, namely it is a global smooth solution on all the space-time.  Existence of such global solutions, certainly exist for $(NS)$. In fact, a constant section $s:M\to W$, is surely a solution for $(NS)$, localized on a equipotential space region. In fact such solution satisfies $\widehat{(NS)}$ iff equations {\em(\ref{constraint-for-constant-section-navier-stokes})} are satisfied.
\begin{equation}\label{constraint-for-constant-section-navier-stokes}
  \left\{
  \begin{array}{l}
    v^k G^j_{jk}=0 \\
    v^k (\partial x_\alpha.G^j_{jk})=0 \\
    v^s R^j_s+\rho(\partial x_i.f)g^{ij}=0\\
    v^kv^p W_{kp}=0.\\
  \end{array}
  \right.
\end{equation}
We have adopted the same symbols used in \cite{PRA1}. Then by using global cartesian coordinates (this is possible for the affine structure of $J^2_4(W)$), we get that $g^{ij}=\delta^{ij}$, $R^j_s=0$ and $W_{kp}=0$. Therefore equations {\em(\ref{constraint-for-constant-section-navier-stokes})} reduce to $\rho(\partial x_k.f)=0$. This means that such constant solutions exist iff they are localized in a equipotential space-region.

Such constant global smooth solutions, even if very simple, can be used to build more complex ones, by using the linearized Navier-Stokes equation at such solutions. Let us denote by $\overline{(NS)}[s]\subset J{\it D}^2(s^*vTW)$ such a linearized PDE at the constant solution $s$. Symilarly to the nonlinear case, we can associate to $\overline{(NS)}[s]$ a linear sub-PDE $\widehat{\overline{(NS)}}[s]\subset\overline{(NS)}[s]$ that is formally integrable and completely integrable. Then in a space-time neighbourhood of a point $q\in \widehat{\overline{(NS)}}[s]$ we can build a smooth solution, say $\nu:M\to s^*vTW$. Since solutions of $\widehat{\overline{(NS)}}[s]$ locally transform solutions of $\widehat{(NS)}$ into other solutions of this last equation, we get that the original constant solution $s$ can be transformed by means of the perturbation $\nu$ into another global solution $s':M\to W$; the perturbation being only localized into a local space-time region. In this way we are able to obtain global space-time smooth solutions $V'\subset (NS)$. (See Fig. \ref{global-nonconstant-smooth-solution-representation}.) Since $V$ and $V'$ are both diffeomorphic to $M$, via the canonical projection $\pi_{2}:J^2_4(W)\to M$, their characteristic classes are all zero: $i_V^*\omega_{i}^{(2)}=i_{V'}^*\omega_{i}^{(2)}=0$, $i\in\{1,2,3,4\}$. Really $H^i(V';\mathbb{Z}_2)=0=H^i(V';\mathbb{Z}_2)$, $\forall i>0$. In the words of Theorem \ref{maslov-index-and-maslov-cycle-relation-solution-pde} we can say that in these global solutions $V$ one has,
\begin{equation}\label{maslov-cycle-ns-compact}
  \Sigma_i(V)_{K}=\varnothing,\,\forall i\in\{1,2,3\}
\end{equation}
for any compact domain $K\subset V$. In {\em(\ref{maslov-cycle-ns-compact})} $\Sigma_i(V)_{K}$ denotes the $i$-Maslov cycle of $V$ inside the compact domain $K\subset V$. (For more details on the existence of such smooth solutions, built by means of perturbations of constant ones, see Appendix A.)\footnote{It is clear that whether we work with the constant solution with zero flow, we get a non-constant
global solutions $V'$ necessarily satisfying the following Clay-Navier-Stokes conditions:

1. $\mathbf{v}(x,t)\in\left[C^\infty(\mathbb{R}^3\times[0,\infty))\right]^3\,,\qquad p(x,t)\in C^\infty(\mathbb{R}^3\times[0,\infty))$

2. There exists a constant $E\in (0,\infty)$ such that $\int_{\mathbb{R}^3} \vert \mathbf{v}(x,t)\vert^2 dx <E$.

For more details on the Navier-Stokes Clay-problem see the following reference: \cite{FEFFERMAN}. Therefore this is another way to prove existence of global smooth solutions when one aims to obtain solutions defined on all the space $\mathbb{R}^3$. Really by varying the localized perturbation one can obtain different initial conditions, and as a by-product global smooth solutions. Such global solutions do not necessitate to be stable at short times, since the symbol of the Navier-Stokes equation is not zero. However by working on the infinity prolongation $\widehat{(NS)}_{+\infty}\subset J^\infty_4(W)$ all smooth solutions can be stabilized at finite times, since for $\widehat{(NS)}_{+\infty}$ one has $(g_2)_{+\infty}=0$. Their average-stability can be studied with the geometric methods given by A. Pr\'astaro, also for global solutions defined on all the space, assuming perturbations with compact support. (See Refs. \cite{PRA6,PRA8,PRA9,PRA10,PRA11,PRA14}.)}
\end{example}

\begin{theorem}[Maslov index for Lagrangian manifolds]\label{maslov-index-for-lagrangian-manifolds}
For $n$-Lagrangian submanifolds of a $2n$-dimensional symplectic manifold $(W,\omega)$, we recognize $i$-Maslov indexes $\beta_i(V)$ and $i$-Maslov cycles $\Sigma_i(V)$, $1\le i\le n-1$. For $i=1$, we recover the Maslov index defined by Arnold.
\end{theorem}
\begin{proof}
After recognized Maslov index for solutions of PDEs (Theorem \ref{maslov-index-and-maslov-cycle-relation-solution-pde}), the first step to follow is to show that Lagrangian submanifolds of $W$ are encoded by a suitable PDE.
Let $\{x^\alpha,y^j\}_{1\le\alpha,\, j\le n}$ be local coordinates in a neighborhood of a point $a\in W$. In this way a $n$-dimensional submanifold $N\subset W$, passing for $a$, can be endowed with local coordinates $\{x^\alpha\}_{1\le\alpha\le n}$. Let us represent $\omega$ in such a coordinate system:
\begin{equation}\label{local-representation-symplectic-form}
  \omega=\sum_{1\le\alpha<\beta\le n}\omega_{\alpha\beta}dx^\alpha\wedge dx^\beta+\sum_{1\le\alpha, j\le n}\bar\omega_{\alpha j}dx^\alpha\wedge dy^j+\sum_{1\le i<j\le n}\hat\omega_{i j}dy^i\wedge dy^j.
\end{equation}
Then the restriction of $\omega$ on a $n$-dimensional submanifold $N\subset W$, with local coordinate $\{x^\alpha\}$, gives the formula (\ref{local-representation-symplectic-form-on-n-dimensional-submanifold}).
\begin{equation}\label{local-representation-symplectic-form-on-n-dimensional-submanifold}
  \omega|_N=\sum_{1\le\alpha<\beta\le n}[\omega_{\alpha\beta}+\sum_{1\le j\le n}(\bar\omega_{\alpha j}y^j_\beta-\bar\omega_{\beta j}y^j_\alpha)+\sum_{1\le i<j\le n}\hat\omega_{i j}(y^i_\alpha y^j_\beta-y^i_\beta y^j_\alpha)]dx^\alpha\wedge dx^\beta.
\end{equation}
Therefore, by imposing that must be $\omega|_N=0$, we see that we can encode $n$-dimensional Lagrangian submanifolds $N\subset $W  by means of solutions, of the PDE $\mathcal{L}_1$ reported in (\ref{pde-lagrangian-n-dimensional-submanifolds}).
\begin{equation}\label{pde-lagrangian-n-dimensional-submanifolds}
 \mathcal{L}_1\subset J^1_n(W)\,:\,\left\{ \begin{array}{l}
  \omega_{rs}(x)+\sum_{1\le j\le n}(\bar\omega_{rj}(x)y^j_s-\bar\omega_{sj}(x)y^j_r)\\
  +\sum_{1\le i<j\le n}(y^i_ry^j_s-y^j_ry^i_s)\hat\omega_{ij}(x)=0\\
  \end{array}\right\}_{1\le r<s\le n}.
\end{equation}
There $\omega_{rs}$ and $\hat\omega_{ij}$ are non-degenerate skew-symmetric $n\times n$ matrices and  $\bar\omega_{rs}$ is a $n\times n$ matrix, all being analytic functions of $\{x^\alpha\}$.\footnote{Let us emphasize that the coefficients $\omega_{rs}$, $\hat\omega_{ij}$ and  $\bar\omega_{rs}$, are related by some first order constraints coming from the condition that $d\omega=0$. However, for the formal integrability of equation (\ref{pde-lagrangian-n-dimensional-submanifolds}) these constraints can be ruled-out.} One can prove that $E_1$ is a formally integrable and completely integrable PDE. In fact one can see that $\pi_{r-1,r}:(E_1)_{+r}\to(E_1)_{+(r-1)}$ are sub-bundles of $\pi_{r+1,r}:J^{r+1}_n(W)\to J^r_n(W)$, $\forall r\ge 1$. Really, for $r=1$ we get

\begin{equation}\label{first-order-prolongation-surjectivity-lagrangian-submanifolds-pde}
 \left\{\begin{array}{l}
      \dim(\mathcal{L}_1)_{+1}=n+n\frac{(+2)(n+1)}{2}-\frac{n(n-1)}{2}-\frac{n^2(n-1)}{2}\\
      \dim(\mathcal{L}_1)=n+n(n+1)-\frac{n(n-1)}{2}\\
       \dim(g_1)_{+1}=\frac{n^2(n+1)}{2}-\frac{n^2(n-1)}{2}\\
 \dim(\mathcal{L}_1)_{+1}= \dim(\mathcal{L}_1)+ \dim(g_1)_{+1}.\\
 \end{array}\right.
\end{equation}
This is enough to state that the short sequence
\begin{equation}\label{short-exact-sequence-first-prolongation-lagrangians-pde}
  \xymatrix{(\mathcal{L}_1)_{+r}\ar[r]^{\pi_{r+1,r}}&(\mathcal{L}_1)_{+(r-1)}\ar[r]&0\\}
\end{equation}
is exact for $r=1$. Since this process can be iterated for any $r>1$, we can state that one can arrive to a certain prolongation where the symbol is involutive. hence the PDE $\mathcal{L}_1$ is formally integrable. Since it is analytic it is also completely integrable.
Then we can apply Theorem \ref{maslov-index-and-maslov-cycle-relation-solution-pde} to $\mathcal{L}_1\subset J^1_n(W)$, to state that there exists Maslov cycles and Maslov indexes for any solution $V\subset(\mathcal{L}_1)_{+1}$, on the first prolongation of $\mathcal{L}_1$. One has the following commutative diagram where all the vertical line are surjectives.
\begin{equation}\label{surjective-mapping-commutative-diagram-lagrangian-submanifolds}
\xymatrix@C=2cm{(\mathcal{L}_1)_{+1}\ar@/_2pc/[dd]_{\pi_{2,0}|_{(\mathcal{L}_1)_{+1}}}\ar[d]^{\pi_{2,1}|_{(\mathcal{L}_1)_{+1}}}\ar@{{(}->}[r]&J^2_n(W)\ar@/^2pc/[dd]^{\pi_{2,0}}\ar[d]_{\pi_{2,1}}\\
  \mathcal{L}_1\ar[d]^{\pi_{1,0}|_{\mathcal{L}_1}}\ar@{{(}->}[r]&J^1_n(W)\ar[d]_{\pi_{1,0}}\\
  W\ar[d]\ar@{=}[r]&W\ar[d]\\
  0&0\\}
\end{equation}

$(\mathcal{L}_1)_{+1}$ is a strong retract of $J^2_n(W)$, hence one has the homotopy equivalence: $$(\mathcal{L}_1)_{+1}\simeq J^2_n(W),$$ that induces isomorphisms on the corresponding cohomology spaces. Therefore, we can recognize $i$-Maslov index còlasses and $i$-Maslov cycle classes on each solution $V\subset (\mathcal{L}_1)_{+1}$.

As a by-product we can apply these results to the symplectic space $(\mathbb{R}^{2n},\omega)$, to recover the same results given by V.I Arnold. (See Example \ref{example-lagrangian-submanifold-a}.) This justifies our Definition \ref{maslov-cycle-solution-pde} and Definition \ref{maslov-index-class-solution-pde} that can be recognized suitable generalizations, in PDE geometry, of analogous definitions given by V. I. Arnold.
\end{proof}

\begin{theorem}[$G$-singular Lagrangian bordism groups]\label{g-singular-lagrangian-bordism-groups}
Let $W$ be a symplectic $2n$-dimesional manifold. Let $G$ be an abelian group. Then the $G$-singular bordism group of $(n-1)$-dimensional compact submanifolds of $W$, bording by means of $n$-dimensional Lagrangian submanifolds of $W$, is given in {\em(\ref{g-singular-lagrangian-bordism-groups-formula})}.
\begin{equation}\label{g-singular-lagrangian-bordism-groups-formula}
  {}^G\Omega^{\mathcal{L}_1}_{\bullet,s}\cong\bar H_\bullet(\mathcal{L}_1;G).
  \end{equation}
$\bullet$\hskip 2pt If ${}^G\Omega^{\mathcal{L}_1}_{\bullet,s}=0$ one has: $\bar Bor_\bullet(\mathcal{L}_1;G)\cong\bar Cyc_\bullet(\mathcal{L}_1;G)$.

$\bullet$\hskip 2pt If $\bar Cyc_\bullet(\mathcal{L}_1;G)$ is a free $G$-module,one has the isomorphism:
$$\bar Bor_\bullet(\mathcal{L}_1;G)\cong {}^G\Omega(\mathcal{L}_1)_{\bullet,s}\bigoplus\bar Cyc_\bullet(\mathcal{L}_1;G).$$
\end{theorem}
\begin{proof}
It is enough to applying Theorem \ref{integral-singular-bordism-groups-of-pde} and Theorem \ref{maslov-index-for-lagrangian-manifolds} to get formula (\ref{g-singular-lagrangian-bordism-groups-formula}).
\end{proof}

\begin{theorem}[Closed weak Lagrangian bordism groups]\label{weak-and-singular-lagrangian-bordism-groups-a}
Let $W$ be a symplectic $2n$-dimesional manifold. Let $G$ be an abelian group. Then the weak $(n-1)$-bordism group of closed compact $(n-1)$-dimensional submanifolds of $W$, bording by means of $n$-dimensional Lagrangian submanifolds of $W$,
is given in {\em(\ref{weak-and-singular-lagrangian-bordism-groups-formula})}.
\begin{equation}\label{weak-and-singular-lagrangian-bordism-groups-formula}
\Omega^{\mathcal{L}_1}_{n-1,w}\cong\bigoplus_{r+s=n-1}H_r(W;{\mathbb Z
}_2)\otimes_{{\mathbb Z}_2}\Omega_s
\cong\Omega^{\mathcal{L}_1}_{n-1}/K^{\mathcal{L}_1}_{n-1,w}\cong
\Omega^{\mathcal{L}_1}_{n-1,s}/K^{\mathcal{L}_1}_{n-1,s,w}.
\end{equation}
 Furthermore, since
$\mathcal{L}_1\subset J^1_n(W)$, has non zero symbols: $g_{1+s}\not=0$,
$s\ge 0$, then $K^{\mathcal{L}_1}_{n-1,s,w}=0$, hence
$\Omega^{\mathcal{L}_1}_{n-1,s}\cong\Omega^{\mathcal{L}_1}_{n-1,w}$.
\end{theorem}

\begin{proof}
From the proof of Theorem \ref{g-singular-lagrangian-bordism-groups} and by using Theorem \ref{formal-integrability-integral-bordism-groups} we get
directly the proof.
\end{proof}

{\bf Warning.}  Lagrangian bordism considered in this paper, namely Theorem \ref{g-singular-lagrangian-bordism-groups} and Theorem \ref{weak-and-singular-lagrangian-bordism-groups-a}, adopts a point of view that is directly related to one where compact (closed) manifolds bording by means of Lagrangian manifolds must be Lagrangian manifolds too. This is, for example the  Lagrangian bordism considered in \cite{BIRAN-CORNEA}. Really these authors work on a manifold $W=\mathbb{R}^2\times M$, where $M$ is a (compact) $2m$-dimensional symplectic manifold $(M,\widehat{\omega})$, and $\mathbb{R}^2$ is endowed with the canonical symplectic form $\omega_{\mathbb{R}^2}=dx\wedge dy$. Thus $W$ is a $2(m+1)$-dimensional symplectic manifold with symplectic form $\omega=\widehat{\omega}\oplus\omega_{\mathbb{R}^2}$. Therefore, one has a natural trivial fiber bundle structure $\pi:W\to \mathbb{R}^2$, with fiber the symplectic manifold $M$. Then one considers bordisms of (closed) compact Lagrangian $m$-dimensional submanifolds of $M$, bording by means of $(m+1)$-dimensional Lagrangian submanifolds of $W$. In such a situation, with respect to the framework considered in Theorem \ref{g-singular-lagrangian-bordism-groups} and Theorem \ref{weak-and-singular-lagrangian-bordism-groups-a} one should specify that $n=m+1$, and that the $n-1=m$ compact manifolds bording with $(n=m+1)$-Lagrangian submanifolds of $W$ must be Lagrangian submanifolds of $M$. In other words the Lagrangian bordism groups considered in \cite{BIRAN-CORNEA} are relative Lagrangian bordism groups, with respect to the submanifold $M\subset W$, in our formulation. However, since $(m+1)$-dimensional Lagrangian submanifolds $V$ of $W$, must necessarily be transverse to the fibers of $\pi:W\to\mathbb{R}^2$, except in the singular points, it follows that the compact (closed) $m$-dimensional manifolds $N_1$ and $N_2$ that they bord, namely $\partial V=N_1\sqcup N_2$, must necessarily be submanifolds of $M$: $N_1,\, N_2\subset M$. Furthermore, by considering that $\omega|_V=0$, it follows that $\widehat{\omega}|_{N_1}=\widehat{\omega}|_{N_2}=0$, hence $N_1$ and $N_2$ must necessarily be Lagrangian submanifolds of $(M,\widehat{\omega})$, as considerd in \cite{BIRAN-CORNEA}. Therefore, our point of view is more general than one adopted in \cite{BIRAN-CORNEA} and recovers this last one when the structure of the symplectic manifold $(W,\omega)$ is of the type $(M\times\mathbb{R}^2,\widehat{\omega}\oplus\omega_{\mathbb{R}^2})$.

\begin{appendices}
\appendix{\bf Appendix A: On global smooth solutions of the Navier-Stokes PDEs}\label{appendixa}
\renewcommand{\theequation}{A.\arabic{equation}}
\setcounter{equation}{0}  

In this appendix we shall explicitly prove a theorem that one has implicitly used in Example \ref{navier-stokes-pdes-and-global-space-time-smooth-solutions}.
\begin{atheorem}\label{theorem-appendix-a}
Any constant smooth solution $s$ of the Navier-Stokes equation $(NS)\subset J\mathcal{D}^2(W)\subset J^2_4(W)$, admits perturbations that identify smooth non-constant solutions of $(NS)\subset J\mathcal{D}^2(W) \subset J^2_4(W)$.
\end{atheorem}

\begin{proof}
We shall use a surgery technique in order to prove this theorem. Let us divide the proof in some lemmas.
\begin{alemma}\label{lemma-a-appendix-a}
Given a smooth constant solution $s$ of $(NS)\subset J\mathcal{D}^2(W)$ we can identify a smooth solution with boundary diffeomorphic to $S^3$ and a compact smooth solution with boundary diffeomorphic to $S^3$ again, such that their canonical projections on $M$ identify an annular domain in $M$.
\end{alemma}
\begin{proof}
Let us consider a compact domain $D\subset M$ identified with a $4$-dimensional disk $D^4$. Set $\partial D^4=S^3$. By fixing a constant solution $s$ of $(NS)\subset J\mathcal{D}^2(W)$, let us denote by $N$ the image into $\widehat{(NS)}$ of $S^3$ by means of $D: N=D^2s(S^3)\subset \widehat{(NS)}\subset J^2_4(W)$. Set $V=D^2s(M)\subset \widehat{(NS)}$ and set
\begin{equation}\label{equation-a-appendix-a}
  \widetilde{V}=(V\setminus D^2s(D^4))\bigcup N\subset\widehat{(NS)}.
\end{equation}

Then $\widetilde{V}$ is a smooth solution of $\widehat{(NS)}$ with boundary $\partial \widetilde{V}=N\cong S^3$.

Let $p_0\in D^4$ be the center of the disk. Since $\widehat{(NS)}\subset J\mathcal{D}^2(W)$ is completely integrable, we can build a smooth (analytic) solution $s_0$ in a neighborhood $U_0\subset D^4$ of $p_0$, such that $\widehat{V}=D^2s_0(U)\subset \widehat{(NS)}$. We can assume that $s_0$ does not coincide with $s$. (Otherwise we could take a different constant value from $s$.) Let us consider in $U_0$ a disk $D^4_0$ centered on $p_0$. Set
\begin{equation}\label{equation-b-appendix-a}
N_0=D^2s_0(\partial D^k_0)\subset \widehat{V}\subset \widehat{(NS)}.
\end{equation}

 Let us consider
 \begin{equation}\label{equation-c-appendix-a}
\widetilde{\widehat{V}}=(\widehat{V}\setminus D^2s_0(D^4_0))\bigcup N_0\subset \widehat{(NS)}.
\end{equation}

Then $\widetilde{\widehat{V}}$ is a smooth solution of $\widehat{(NS)}$ with boundary $\partial\widetilde{\widehat{V}}=N_0\cong S^3$.
Of course the projections of $N$ and $N_0$ on $M$ via the canonical projection $\pi_{2}:J\mathcal{D}^2(W)\to M$, identify an annular domain in $M$.
\end{proof}
\begin{alemma}\label{lemma-a-appendix-b}
The solutions $\widetilde{V}$ and $\widetilde{\widehat{V}}$ considered in Lemma \ref{lemma-a-appendix-a} identify a connected smooth solution of $(NS)\subset J^2_4(W)$.
\end{alemma}
\begin{proof}
Since both solutions $\widetilde{V}$ and $\widetilde{\widehat{V}}$ are smooth solutions we can consider their $\infty$-prolongations and look to them inside $\widehat{(NS)}_{+\infty}$. Now their boundary are both diffeomorphic to $S^3$, therefore must exist a smooth solution $\mathop{V}\limits^{\infty}\subset \widehat{(NS)}_{+\infty}\subset J^\infty_4(W)$ such that $\partial \mathop{V}\limits^{\infty}=N_{+\infty}\bigcup(N_0)_{+\infty}$. In fact, from the commutative diagram \ref{commutative-diagram-relating-bordism-groups-in-pde} and Theorem \ref{formal-integrability-integral-bordism-groups} we get the exact commutative diagram (\ref{bordism-sequence-smooth-solutions-navier-stokes-pde}).

\begin{equation}\label{bordism-sequence-smooth-solutions-navier-stokes-pde}
  \xymatrix{0\ar[dr]&0&&&\\
  &\overline{K}^{\widehat{(NS)}}_3\ar[u]\ar[dr]&&0\ar[d]&\\
  0\ar[r]&K^{\widehat{(NS)}}_{3,s}\ar[u]\ar[r]&\Omega_3^{\widehat{(NS)}}\ar[rd]\ar[r]&\Omega_{3,s}^{\widehat{(NS)}}\ar[d]\ar[r]&0\\
  &0\ar[u]\ar@{=}[rr]&&\Omega_3\ar[rd]\ar[d]&\\
  &&&0&0\\}
\end{equation}
where $\overline{K}^{\widehat{(NS)}}_3\cong K^{\widehat{(NS)}}_{3,s}$ distinguishes between non-diffeomorphic closed $3$-dimensional integral smooth submanifolds of $\widehat{(NS)}$. In fact $\Omega_{3,s}^{\widehat{(NS)}}\cong\Omega_3=0$.\footnote{This is related to the fact that the Navier-Stokes equation is an {\em extended $0$-crystal PDE}. (See \cite{PRA8,PRA9,PRA10,PRA11,PRA14,PRA16,PRA017}.) In Tab. \ref{some-unoriented-smooth-bordism-groups} are reported some un-oriented smooth bordism groups $\Omega_n$, $0\le n\le 3$, useful to calculate $\Omega_{3,s}^{\widehat{(NS)}}$, according to Theorem \ref{formal-integrability-integral-bordism-groups}.} Since the Cartan distribution $\mathbf{E}_{\infty}\subset T\widehat{(NS)}_{+\infty}$ is $4$-dimensional, it follows that $\mathop{V}\limits^{\infty}$ smoothly solders with the solutions $\mathop{V}\limits^{\infty}$, $\widetilde{V}$ and $\widetilde{\widehat{V}}$. In this way
\begin{equation}\label{equation-d-appendix-a}
  X=\widetilde{V}_{+\infty}\bigcup_{N_{+\infty}}\mathop{V}\limits^{\infty}\bigcup_{(N_0)_{+\infty}}\widetilde{\widehat{V}}_{+\infty}
\end{equation}

is a smooth solution of $\widehat{(NS)}_{+\infty}$, hence of $(NS)$.
\end{proof}
To conclude the proof let us assume that $\mathop{V}\limits^{\infty}$ can be realized by means of a section $s_{\infty}$ of $\pi:W\to M$, namely $\mathop{V}\limits^{\infty}=D^\infty (s_{\infty})(A)$, where $A$ is the annular domain above considered. Thus we can say that the solution $X=D^{\infty}\bar s(M)$ for some smooth global section $\bar s$ of $\pi:W\to M$.  Then taking into account of the affine structure of $W$ we can state that $\bar s=s+\nu$, where $\nu$ is a smooth perturbation of $s$ on the disk $D^4$, such that $\nu|_{S^3}=0$ and

\begin{equation}\label{equation-e-appendix-a}
 \mathop{\lim}\limits_{p\to S^3\, \hbox{\rm{\footnotesize(from inside)}}}\nu(p)=0,\, \nu|_{\complement D^4}=0.
 \end{equation}

In other words the perturbation is of the type pictured in Fig. \ref{global-nonconstant-smooth-solution-representation}.

So we have proved the following lemma.
\begin{alemma}\label{lemma-b-appendix-b}
When perturbations of $\widehat{(NS)}$ are realized by means of smooth solutions of the corresponding linearized Navier-Stokes PDE, $\widehat{\overline{(NS)}}[s]\subset J\mathcal{D}^2(s^*vTW)$, the completely integrable part of $\overline{(NS)}[s]$, then the identified solutions of $(NS)\subset J^2_4(W)$ are also smooth solutions of $(NS)\subset J\mathcal{D}^2(W)\subset J^2_4(W)$, namely they are identified with smooth sections of $\pi:W\to M$.
\end{alemma}

\begin{table}[t]
\caption{Un-oriented smooth bordism groups $\Omega_n$, $0\le n\le 3$.}
\label{some-unoriented-smooth-bordism-groups}
\begin{tabular}{|l|l|}
\hline
\hfil{\rm{\footnotesize $n$}}\hfil&\hfil{\rm{\footnotesize $\Omega_n$}}\hfil\\
\hline\hline
\hfil{\rm{\footnotesize $0$}}\hfil&\hfil{\rm{\footnotesize $\mathbb{Z}_2$}}\hfil\\
\hline
\hfil{\rm{\footnotesize $1$}}\hfil&\hfil{\rm{\footnotesize $0$}}\hfil\\
\hline
\hfil{\rm{\footnotesize $2$}}\hfil&\hfil{\rm{\footnotesize $\mathbb{Z}_2$}}\hfil\\
\hline
\hfil{\rm{\footnotesize $3$}}\hfil&\hfil{\rm{\footnotesize $0$}}\hfil\\
\hline
\end{tabular}
\end{table}

Whether, instead, $\mathop{V}\limits^{\infty}$ is a smooth solution that cannot be globally represented by means of a section of $\pi:W\to M$, then it means that there are in $\mathop{V}\limits^{\infty}$, and hence in its projection into $W$, some pieces that climb on the fibers of $\pi:W\to M$. In such a case we can continue to state that $X$ is obtained  by a perturbation of $V$ inside the compact domain $D$, but the perturbation is a singular solution of the linearized Navier-Stokes PDE at the constant section $s$. Therefore in such a case it should not be possible represent $X$ as a smooth section of $\pi:W\to M$. This shows the necessity to realize the perturbation of $s$ by means of a smooth solution $\nu$ of the linearized equation $\widehat{\overline{(NS)}}[s]\subset J\mathcal{D}^2(s^*vTW)$, such that $\nu|_{S^3}=0$ and $\nu|_{\complement D^4}=0$, in order the perturbed solution $X$ should be identified with a global non-constant section of $\pi:W\to M$.

\begin{figure}[t]
\centerline{\includegraphics[width=8cm]{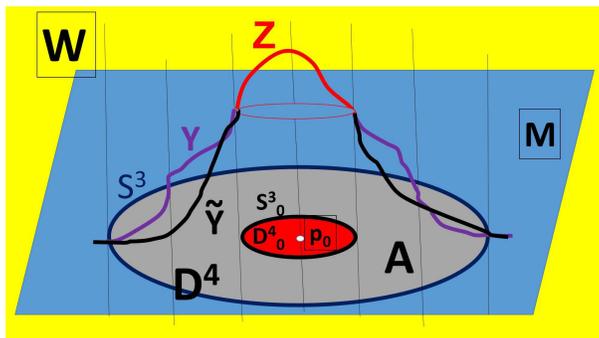}}
\renewcommand{\figurename}{Fig.}
\caption{Deformation of a smooth solution ($Y$) of $\widehat{(NS)}\subset J^2_4(W)$, climbing along the fibers of $\pi_2:\widehat{(NS)}\to M$, into a smooth solution ($\widetilde{Y}$) of $\widehat{(NS)}\subset J\mathcal{D}^2(W)\subset J^2_4(W)$. This is possible since the Navier-Stokes PDE is an affine fiber bundle over $(C)$, and its symbol is not zero: $\dim (g_2)_q=46$, $\forall q\in(NS)$, $\dim (\widehat{g_2})_q=42$, $\forall q\in\widehat{(NS)}$. In the picture $Z=\pi_{2,0}(\widehat{V})\subset W$.}
\label{global-nonconstant-smooth-solution-representation-a}
\end{figure}

On the other hand since $\widetilde{V}_{+\infty}$ and $\widetilde{\widehat{V}}_{+\infty}$ are both regular solutions with respect to the canonical projection $\pi_{\infty}:J\mathcal{D}^\infty(W)\to M$, and $(NS)\subset J\mathcal{D}^2(W)$ is an affine fiber bundle over its projection at the first order, with non-zero symbol, it follows that we can deform any eventual piece climbing on the fibers in such a way to obtain a regular solution with respect to the projection $\pi:W\to M$. Therefore, the projection $Y\subset\widehat{ (NS)}$ of $\mathop{V}\limits^{\infty}$ into $\widehat{ (NS)}$, can be eventually deformed into a regular solution, $\widetilde{Y}$, with respect to the projection $ \pi_{2}$. (See Fig. \ref{global-nonconstant-smooth-solution-representation-a}.) In this way the projection of $\widetilde{Y}$ into $W$, smoothly relates regular smooth submanifolds that project on two domains of $M$ that are outside the annular domain $A$, but that are disconnected each other. Thus $\widetilde{Y}$ identifies a smooth $4$-dimensional manifold transverse to the fibers of $\pi:W\to M$. By conclusion $\widetilde{Y}^{+\infty}$ is necessarily a regular solution of $\widehat{(NS)}_{+\infty}\subset J\mathcal{D}^\infty(W)$. Therefore it can be obtained by a perturbation of the constant solution $s$, by means of a smooth solution of $\widehat{\overline{(NS)}}[s]\subset J\mathcal{D}^2(s^*vTW)$. (See Fig. \ref{global-nonconstant-smooth-solution-representation}.)
\end{proof}

\appendix{\bf Appendix B: On the Legendrian bordism}\label{appendixb}
\renewcommand{\theequation}{B.\arabic{equation}}
\setcounter{equation}{0}  

Similarly to the way we considered Lagrangian bordism in this paper, we can also formulate Legendrian bordism. Let us in this appendix recall some basic definitions and sketch  only some steps on.  Really on a $(2n+1)$-dimensional manifold $W$, endowed with a {\em contact structure}, namely a $1$-differential form $\chi$, such that $\chi(p)\wedge d\chi(p)^n\not=0$, $\forall p\in W$, there exists a {\em characteristic vector field} $v:W\to TW$, i.e., the generator of the $1$-dimensional annihilator of $d\chi$: $v\rfloor d\chi=0$ and $v\rfloor\chi =1$. Furthermore on $W$ there exists also a {\em contact distribution}, namely a $2n$-dimensional distribution $\mathbf{B}=\bigcup_{p\in W}\mathbf{B}_p$, $\mathbf{B}_p=\ker(\chi(p))\subset T_pW$. One has the following properties.
\begin{bproposition}\label{properties-contact-distribution-appendix-b}
The following propositions hold.

{\em(bi)} $d\chi(p)|_{\mathbf{B}_p}$, $\forall p\in W$, is {\em nondegenerate}, i.e., if $d\chi(\zeta,\xi)=0$, $\forall \zeta\in\mathbf{B}_p$, and $\forall \xi\in \mathbf{B}_p$, then $\zeta=0$.

 {\em(bii)} $TW=\mathbf{B}\bigoplus <v>$.

 {\em(biii) (Darboux's theorem)}  $\mathbf{B}\to W$ is a symplectic vector bundle with symplectic form $d\chi|_{\mathbf{B}}$.

 {\em(biv)}  With respect to local coordinates $\{x^\alpha, y_\alpha, z\}$ on $W$, $\chi $ assumes the following form:\footnote{All contact strcture forms on $W$ are locally diffeomorphic.}
 \begin{equation}\label{canonical-local-form-contact-structure-appendix-b}
   \chi= dz-y_\alpha dx^\alpha.
 \end{equation}
\end{bproposition}

\begin{bproposition}\label{properties-integral-manifolds-contact-structure}
{\em Integral manifold} of a contact structure $(W,\chi)$, is a submanifold $N\subset W$, such that $\chi|_{N}=0$, (or equivalently $T_pN\subset \mathbf{B}_p$, $\forall p\in N$). One has
\begin{equation}\label{dimension-integral-manifold-contact-structure-appendix-b}
  \dim N<\frac{1}{2}(2n+1).
\end{equation}
$\bullet$\hskip 2pt {\em Legendrian submanifolds} of $(W,\chi)$ are integral submanifolds $N$ of maximal dimension: $\dim N=n$.
\end{bproposition}

\begin{bdefinition}\label{legendrian-bundles-appendix-b}
A {\em Legendrian bundle} $\pi:W\to M$,  is a fiber bundle  with $\dim W=2n+1$, $\dim M=n+1$, and endowed with a contact structure $(W,\chi)$, such that each fiber $W_p$ is a Legendrian submanifold, namely $\chi|_{W_p}=0$, $\forall p\in M$.

$\bullet$\hskip 2pt If $L\subset W$ is a Legendrian submanifold of $W$, ($\chi|L=0$, $\dim L=n$), its {\em front} is $\pi(L)=X\subset M$. Singularities of $\pi|_L:L\to M$ are called {\em Legendrian singularities}. The front $X$ of a Legendrian submanifold is a $n$-dimensional submanifold of $M$, with eventual singularities.
\end{bdefinition}
Similarly to the Lagrangian submanifolds of symplectic manifolds, we can characterize Legendrian submanifolds of a contact manifold by means of suitable PDEs. In fact we have the following.

\begin{btheorem}\label{legendrian-pde-appendix-b}
Given a contact structure on a $(2n+1)$-dimensional manifold $(W,\chi)$, its Legendrian submanifolds are solutions of a first order, involutive, formally integrable and completely integrable PDE.

$i$-Maslov indexes and $i$-Maslov cycles, $1\le i\le n-1$, can be recognized for such solutions.
\end{btheorem}

\begin{proof}

Let $\{x^\alpha,y_\alpha,z\}_{1\le\alpha\le n}$ be local coordinates on $W$. Then Legendrian submanifolds of $W$ are the $n$-dimensional submanifolds of $W$ that satisfy the PDE reported in (\ref{legendrian-pde-appendix-b-ocal-expression}).

\begin{equation}\label{legendrian-pde-appendix-b-ocal-expression}
\mathcal{L}eg\subset J^1_n(W):\, \{z_\beta-y_\beta=0\}
\end{equation}
where $\{x^\alpha,y_\alpha,z,y_{\alpha\beta},z_\beta\}_{1\le\alpha,\beta\le n}$ are local coordinates on $J^1_n(W)$.
The first prolongation of $\mathcal{L}eg$ is given in (\ref{first-prolongation-legendrian-pde-appendix-b-ocal-expression}).\footnote{Let us note that the equations of second order in (\ref{first-prolongation-legendrian-pde-appendix-b-ocal-expression}) are not all linearly independent. In fact, by considering that must be $z_{\beta\gamma}=z_{\gamma\beta}$, we get, by difference of the equations $z_{\beta\gamma}-y_{\beta\gamma}=0$ and $z_{\gamma\beta}-y_{\gamma\beta}=0$: $-y_{\beta\gamma}+y_{\gamma\beta}=0$, namely it must hold the symmetry under the exchange of indexes in in $y_{\gamma\beta}$. Therefore the number of independent equations for $\mathcal{L}eg_{+1}$ is $(n+\frac{n(n+1)}{2})$.}
\begin{equation}\label{first-prolongation-legendrian-pde-appendix-b-ocal-expression}
\mathcal{L}eg_{+1}\subset J^2_n(W):\left\{
\begin{array}{l}
 z_\beta-y_\beta=0\\
  z_{\beta\gamma}-y_{\beta\gamma}=0
\end{array}\right\}.
\end{equation}

Then one can see that
\begin{equation}\label{dimension-surjectivity-appendix-b}
\left\{\begin{array}{ll}
[\dim(\mathcal{L}eg_{+1})&=\frac{n+1}{2}(n^2+2n+2)]=[\dim(\mathcal{L}eg)=(n+1)^2]\\
&+[\dim((g_1)_{+1})=\frac{(n+1)n^2}{2}].\\
\end{array}\right.
\end{equation}
Therefore, one has the surjectivity $\mathcal{L}eg_{+1}\to \mathcal{L}eg$. Furthermore, one can see that the symbol $g_1$ is involutive. In fact one has
\begin{equation}\label{involutivness-symbols-legendrian-pde}
  \left\{
  \begin{array}{ll}
  [\dim(((g_1)_{+1}))=\frac{n^2(n+1)}{2}&=[\dim(g_1)=n^2]+[\dim(g_1^{(1)})=n^2-n)]\\
  &+[\dim(g_1^{(2)})=n^2-2n)]\\
  &+\cdots+[\dim(g_1^{(n-1)})=n^2-n(n-1)]\\
  &=\frac{n^2(n+1)}{2}.\\
  \end{array}
  \right.
\end{equation}

We have used the formula $1+2+3+\cdots+(n-1)=\frac{n(n-1)}{2}$. This is enough to state that $\mathcal{L}eg$ is formally integrable and being analytic it is also completely integrable.

Let us also remark that $\mathcal{L}eg$ is a strong retract of $J^1_n(W)$, thereore one has the homotopic equivalence $J^1_n(W)\simeq\mathcal{L}eg$ that induces isomorphisms between the corresponding cohomology groups. Then by using Theorem \ref{maslov-index-and-maslov-cycle-relation-solution-pde} we can state that on each solution of $\mathcal{L}eg$ we are able to recognize $i$-Maslov indexes and $i$-Maslov cycles.
\end{proof}
\begin{bdefinition}\label{legendrian-bordism-appendix-b}
Let $W$ be a $(2n+1)$-dimensional contact manifold $(W,\chi)$. A Legendrian bordism is a $n$-dimensional Legendrian submanifold bording compact $(n-1)$-dimensional integral submanifolds of $W$.
\end{bdefinition}

\begin{bexample}\label{example-contact-manifold-appendix-b}
Let $M$ be a $n$-dimensional manifold. The derivative space $$J\mathcal{D}(M,\mathbb{R})\cong T^*M\times \mathbb{R}$$ has a canonical contact form $\chi=dy-y_\alpha dx^\alpha$, where $x^\alpha$ are local coordinates on $M$ and $y$ is a coordinate on $\mathbb{R}$. This is just the Cartan form on the derivative space $J\mathcal{D}(E)$, $E=M\times \mathbb{R}$, with respect to the fibration $\pi:E\to M$. The corresponding contact distribution coincides with the Cartan distribution $\mathbf{E}_1(E)\subset TJ\mathcal{D}(E)$. Every solution is a Legendrian submanifold.  Therefore in such a case Legendrian bordism are identified with solutions bording $(n-1)$-dimensional integral submanifolds.
\end{bexample}
\begin{bremark}\label{final-remark-appendix-b}
From above results we can directly reproduce results similar to Theorem \ref{g-singular-lagrangian-bordism-groups} and Theorem \ref{weak-and-singular-lagrangian-bordism-groups-a} also for singular Legendrian bordism groups. More precisely one has the exact commutative diagram reported in (\ref{exact-commutative-diagram-homotopy-equivalence-legendrian-pde}, where the top horizontal line is an homotopy equivalence.
\begin{equation}\label{exact-commutative-diagram-homotopy-equivalence-legendrian-pde}
  \xymatrix@C=2cm{\mathcal{L}eg\ar@{-}[r]^{\sim}\ar[d]&J^1_n(W)\ar[d]\\
  W\ar[d]\ar@{=}[r]&W\ar[d]\\
  0&0\\}
\end{equation}
We get the following isomorphisms:
\begin{equation}\label{isomorphisms-legendria-integral-planes}
  \left\{
  \begin{array}{ll}
  H^1(I(\mathcal{L}eg);\mathbb{Z}_2)&\cong H^1(W;\mathbb{Z}_2)\bigoplus \mathbb{Z}_2[\omega_1^{(1)}]\\
    H^i(I(\mathcal{L}eg);\mathbb{Z}_2)&\cong H^i(W;\mathbb{Z}_2)\\
    &\bigoplus_{1\le p\le i-1}H^{i-p}(W;\mathbb{Z}_2)\bigotimes_{\mathbb{Z}_2}H^p(F_1;\mathbb{Z}_2)\\
    &\bigoplus \mathbb{Z}_2[\omega_1^{(1)},\cdots,\omega_i^{(1)}].\\
  \end{array}
  \right.
\end{equation}
Then the map $i_V:V\to I(\mathcal{L}eg)$ induces the following morphism
\begin{equation}\label{morphism-on-cohomology-legendria-integral-planes}
(i_V)_*:\mathbb{Z}_2[\omega_1^{(1)},\cdots,\omega_i^{(1)}]\to H^i(V;\mathbb{Z}_2),\, 1\le i\le n-1.
\end{equation}
Set $\beta_i(V)=(i_V)_*(\omega_i^{(1)})$ that is the $i$-Maslov index of the Legendrian manifold $V$. We get $\beta_i(V)\bigcap[V]=[\Sigma_i(V)]$ that relates the $i$-Maslov index of $V$ with its $i$-Maslov cycle.
\end{bremark}

\begin{btheorem}[$G$-singular Legendrian bordism groups]\label{g-singular-legendrian-bordism-groups}
Let $W$ be a contact $(2n+1)$-dimesional manifold. Let $G$ be an abelian group. Then the $G$-singular bordism group of $(n-1)$-dimensional compact submanifolds of $W$, bording by means of $n$-dimensional Legendrian submanifolds of $W$, is given in {\em(\ref{g-singular-legendrian-bordism-groups-formula})}.
\begin{equation}\label{g-singular-legendrian-bordism-groups-formula}
  {}^G\Omega^{\mathcal{L}eg}_{\bullet,s}\cong\bar H_\bullet(\mathcal{L}eg;G).
  \end{equation}
$\bullet$\hskip 2pt If ${}^G\Omega^{\mathcal{L}eg}_{\bullet,s}=0$ one has: $\bar Bor_\bullet(\mathcal{L}eg;G)\cong\bar Cyc_\bullet(\mathcal{L}eg;G)$.

$\bullet$\hskip 2pt If $\bar Cyc_\bullet(\mathcal{L}eg;G)$ is a free $G$-module,one has the isomorphism:
$$\bar Bor_\bullet(\mathcal{L}eg;G)\cong {}^G\Omega(\mathcal{L}eg)_{\bullet,s}\bigoplus\bar Cyc_\bullet(\mathcal{L}eg;G).$$
\end{btheorem}

\begin{btheorem}[Closed weak Legendrian bordism groups]\label{weak-and-singular-legendrian-bordism-groups-a}
Let $W$ be a contact $(2n+1)$-dimesional manifold. Let $G$ be an abelian group. Then the weak $(n-1)$-bordism group of closed compact $(n-1)$-dimensional submanifolds of $W$, bording by means of $n$-dimensional Legendrian submanifolds of $W$,
is given in {\em(\ref{weak-and-singular-legendrian-bordism-groups-formula})}.
\begin{equation}\label{weak-and-singular-legendrian-bordism-groups-formula}
\Omega^{\mathcal{L}eg}_{n-1,w}\cong\bigoplus_{r+s=n-1}H_r(W;{\mathbb Z
}_2)\otimes_{{\mathbb Z}_2}\Omega_s
\cong\Omega^{\mathcal{L}eg}_{n-1}/K^{\mathcal{L}eg}_{n-1,w}\cong
\Omega^{\mathcal{L}eg}_{n-1,s}/K^{\mathcal{L}eg}_{n-1,s,w}.
\end{equation}
 Furthermore, since
$\mathcal{L}eg\subset J^1_n(W)$ has non zero symbols then $K^{\mathcal{L}eg}_{n-1,s,w}=0$, hence
$\Omega^{\mathcal{L}eg}_{n-1,s}\cong\Omega^{\mathcal{L}eg}_{n-1,w}$.
\end{btheorem}

\end{appendices}

\end{document}